\documentclass[12pt,mathscr,a4paper,oneside]{amsart}
\usepackage[a4paper]{geometry}
\usepackage{amssymb,eucal}
\usepackage[vcentermath,noautoscale]{youngtab}
\usepackage[UglyObsolete]{diagrams} 

\usepackage{mathptmx} 
\DeclareMathAlphabet{\mathcal}{OMS}{cmsy}{m}{n} 

\theoremstyle{plain}
\newtheorem{theorem}{Theorem}[section]
\newtheorem{proposition}[theorem]{Proposition}
\newtheorem{lemma}[theorem]{Lemma}
\newtheorem{corollary}[theorem]{Corollary}

\theoremstyle{remark}
\newtheorem*{remark}{Remark}
\newtheorem*{remarks}{Remarks}
\newtheorem*{definitions}{Definitions}

\newcommand{\set}[2]{\ensuremath{\{ #1 \>|\> #2 \}}}

\begin{document}

\title[Low-dimensional cohomology and analogs of the Riemann tensor]
{Low-dimensional cohomology of current Lie algebras and
analogs of the Riemann tensor for loop manifolds}
\author{Pasha Zusmanovich}
\address{}
\email{pasha.zusmanovich@gmail.com}
\date{February 26, 2003; last minor revision January 31, 2018}
\thanks{arXiv:math/0302334; Lin. Algebra Appl. \textbf{407} (2005), 71--104}
\subjclass[2000]{Primary 17B56, Secondary 15A69, 17B20, 17B70, 32M15, 53C10, 58H15}

\begin{abstract}
We obtain formulas for the first and second cohomology groups of 
a general current Lie algebra with coefficients in the ``current'' module, 
and apply them to compute structure functions for
manifolds of loops with values in compact Hermitian symmetric spaces.
\end{abstract}

\maketitle

\section*{Introduction}

We deal with the low-dimensional cohomology of current Lie algebras with 
coefficients in the ``current module''.
Namely, let $L$ be a Lie algebra, $M$ an $L$-module, $A$ an
associative commutative algebra with unit, $V$ a symmetric unital
$A$-module. Then the Lie algebra structure on $L\otimes A$ and the
$L\otimes A$-module structure on $M\otimes V$ are defined via obvious formulas:
\begin{align*}
[x \otimes a, y\otimes b] &= [x,y] \otimes ab,                      \\
(x\otimes a) \bullet (m\otimes v) &= (x\bullet m) \otimes (a\bullet v)
\end{align*}
for any $x,y \in L$, $m \in M$, $a,b \in A$, $v \in V$, where $\bullet$ denotes,
by abuse of notation, a respective module action.

The aim of this paper is twofold.
First, we want to demonstrate that the problem
of description of such cohomology in terms of the tensor factors $L$ and $A$
probably does not have an adequate general solution, 
as even a partial answer for the two-dimensional cohomology seems to be 
overwhelmingly complex.
Second, we want to demonstrate, nevertheless, computability of this cohomology
in some cases and its application to some differential geometric questions.

In \S 1 we establish an elementary result from linear algebra which will be useful
in the course of subsequent algebraic manipulations. 
In \S 2 we get a formula for the first cohomology group. 
In \S 3 we compute the second cohomology group in two cases -- where $L$ is abelian
and where $L$ acts trivially on the whole cohomology
group $H^2(L \otimes A, M \otimes V)$. 
At the end of this section, we present a list of 13 types of
2-cocycles (so-called cocycles of rank 1, generated by decomposable elements in
the tensor product) in the general case. However, this list is a priori not complete.
In \S 4 a certain spectral sequence is sketched, which may provide a more
conceptual framework for computations in preceding sections. However, we do not go
into details and other sections are not dependent on that one.
The last \S 5 is devoted to an application. We show how to derive from previous computations
obstructions to integrability (structure functions)
of certain canonical connections on the manifolds of loops with values in 
compact Hermitian symmetric spaces.

One should note that the result about the first cohomology group 
(in particular, about derivations of the current Lie algebra) can be found
in different forms in the literature and is a sort of folklore, 
and partial results on the second cohomology were obtained by 
Cathelineau \cite{cat}, Haddi \cite{haddi-3}, Lecomte and Roger \cite{roger} 
and the author \cite{deformations}.
However, all these results do not provide the whole generality we need, as various
restrictions, notably the zero characteristic of the 
ground field and perfectness of the Lie algebra $L$ were imposed. 
Moreover, as we see in \S 5, the case in a sense opposite to the case of perfect $L$, 
namely, the case of abelian $L$, does lead to some interesting application
(first considered by Poletaeva).

The technique used is highly computational and linear-algebraic in nature and based on 
applying various symmetrization operators to the cocycle equation.

\section*{Notations}

The ground field $K$ is assumed to be arbitrary field of characteristic $\ne 2,3$ 
in \S 1--4, and $\mathbb C$ in \S 5.

$H^n(L,M)$, $C^n(L,M)$, $Z^n(L,M)$, $B^n(L,M)$ stand, respectively, for the spaces of 
cohomology, cochains, cocycles and coboundaries of a Lie algebra $L$
with coefficients in a module $M$.

$M^L = \set{m\in M}{x\bullet m =0 \>\text{for any}\> x\in L}$ is a submodule of 
$L$-invariants.

If $M,N$ are two $L$-modules, $Hom (M,N)$ bears a standard $L$-module structure via
($x\bullet \varphi)(m) = \varphi(x\bullet m) - x\bullet \varphi(m)$ for $x\in L, m\in M$,
and $Hom_L (M,N)$ is another notation for $Hom(M,N)^L$.

$S^n(A,V)$ stands for the space of $n$-linear maps 
$A \times \dots \times A \to V$, symmetric in all arguments.

$\wedge^n(V)$ and $T^n(V)$ stand, respectively, for the spaces of 
$n$-fold skew and tensor products of a module $V$.

$Har^n(A,V)$ and $\mathscr Z^n(A,V)$ stand, respectively, for the spaces of
Harrison cohomology and Harrison cocycles of an associative commutative algebra $A$ with coefficients 
in a module $V$ (for $n=2$, these are just symmetric Hochschild cocycles; 
see \cite{harrison}, where this cohomology was introduced, and \cite{GS} for a more modern
treatment).

\begin{sloppypar}
$Der(A)$ denotes the derivation algebra of an algebra $A$. More generally, 
$Der(A,V)$ denotes the space of derivations of $A$ with values in a $A$-module $V$.
\end{sloppypar}

All other (nonstandard and unavoidably numerous) notations for different spaces of 
multilinear mappings and modules are defined as they introduced in the text.

The symbol $\curvearrowright$ after an expression refers to the sum of all cyclic permutations
(under $S(3)$) of letters and indices occurring in that expression.

\section{A lemma from linear algebra}

If either both $L$ and $M$ or both $A$ and $V$ are finite-dimensional, 
then each cocycle $\Phi$ in $Z^n(L\otimes A, M\otimes V)$ 
can be represented as an element of 
$Hom\,(L^{\otimes n}, M) \otimes Hom\,(A^{\otimes n}, V)$:
\begin{equation}\tag{1.1}
\Phi = \sum_{i\in I} \varphi_i \otimes \alpha_i
\end{equation}
where $\varphi_i, \alpha_i$ are $n$-linear mappings 
$L \times \dots \times L \to M$ and $A \times \dots \times A \to V$
respectively. We restrict our considerations to this case. The minimal possible
number $|I|$ such that the cocycle $\Phi$ can be written in the form (1.1) 
will be called the \textit{rank of cocycle}.

Representing $H^n(L \otimes A, M \otimes V)$
in terms of pairs $(L, M)$ and $(A, V)$, we encounter conditions such as
\begin{equation}\tag{1.2}
\sum_{i\in I} S\varphi_i \otimes T\alpha_i = 0,
\end{equation}
where $S$ and $T$ are some linear operators defined on the spaces of $n$-linear
mappings $L \times \dots \times L \to M$ and $A \times \dots \times A \to V$,
respectively.

For example, the substitution $a_1 = \dots = a_{n+1} = 1$ in the cocycle equation
$d\Phi(x_1\otimes a_1, \dots, x_{n+1}\otimes a_{n+1}) = 0$, where 
$\Phi$ is as in (1.1), yields
$$\notag
\sum_{i\in I} d\varphi_i(x_1, \dots, x_{n+1}) \otimes \alpha_i(1, \dots, 1) = 0.
$$

Another example: applying the symmetrization operator $Y$
with respect to the letters $x_1, \dots, x_{n+1}$, to the cocycle equation, we get:
$$\notag
\sum_{i\in I} \left( Y(x_1 \bullet \varphi_i(x_2, \dots, x_{n+1})) \otimes
                     \sum_{j=1}^{n+1} (-1)^j a_j \bullet 
                     \alpha_i(a_1, \dots, \widehat{a_j}, \dots, a_{n+1}) 
              \right) = 0.
$$

So, suppose that a condition of type (1.2) holds. Since
$$ 
Ker (S\otimes T) = Hom\,(L^{\otimes n}, M) \otimes Ker\,T   +
                   Ker\,S \otimes Hom\,(A^{\otimes n}, V),
$$
it follows that replacing $\alpha_i$'s and $\varphi_i$'s by appropriate linear
combinations, one can find a decomposition of the set of indices $I = I_1 \cup I_2$
such that
\begin{equation}\tag{1.3}
S\varphi_i = 0, \quad i \in I_1 
\qquad\text{and}\qquad
T\alpha_i = 0, \quad i \in I_2.
\end{equation}

Suppose that another equality of type (1.2) holds:
\begin{equation}\tag*{$(1.2^\prime)$}
\sum_{i\in I} S^\prime\varphi_i \otimes T^\prime\alpha_i = 0.
\end{equation}
Then it determines a new decomposition $I = I_1^\prime \cup I_2^\prime$
such that 
$S^\prime \varphi_i = 0$ if $i \in I_1^\prime$ and
$T^\prime \alpha_i = 0$ if $i \in I_2^\prime$.
It turns out that it is possible to replace $\varphi_i$'s and $\alpha_i$'s by their linear
combinations so that both decompositions will hold simultaneously.

\begin{lemma}\label{1.1}
\begin{sloppypar}
Let $U, W$ be two vector spaces, 
$S, S^\prime \in Hom\,(U, \>\cdot\>)$, $T, T^\prime \in Hom\,(W, \>\cdot\>)$. Then
\begin{multline}\notag
Ker\,(S \otimes T) \cap Ker\,(S^\prime \otimes T^\prime) \\ \simeq
(Ker\,S \cap Ker\,S^\prime) \otimes W + 
Ker\,S                      \otimes Ker\,T^\prime   + 
Ker\,S^\prime               \otimes Ker\,T          +
U                           \otimes (Ker\,T \cap Ker\,T^\prime).
\end{multline}
\end{sloppypar}
\end{lemma}

\begin{proof}
Since $Ker(S \otimes T) = Ker\,S \otimes W + U \otimes Ker\,T$ and analogously for
$Ker (S^\prime \otimes T^\prime)$, the equality to prove is a particular case of
\begin{multline}\tag{1.4}
(U_1 \otimes W + U \otimes W_1) \cap 
(U_2 \otimes W + U \otimes W_2) \\ = 
(U_1 \cap U_2) \otimes W + 
U_1 \otimes W_2 + 
U_2 \otimes W_1 + U \otimes (W_1 \cap W_2)
\end{multline}
provided $U_1, U_2$ and $W_1, W_2$ are subspaces of $U$ and $W$ respectively.

Assume for the moment that $U_1 \cap U_2 = W_1 \cap W_2 = 0$. Then expressing
$U = U_1 \oplus U_2 \oplus U^\prime$ and 
$W = W_1 \oplus W_2 \oplus W^\prime$ for some subspaces $U^\prime, W^\prime$
and substituting this in the left side of (1.4), we get:
\begin{multline}\notag
(U_1      \otimes W    \>\oplus\> 
 U_2      \otimes W_1  \>\oplus\>
 U^\prime \otimes W_1)        \cap 
(U_1      \otimes W_2  \>\oplus\>
 U_2      \otimes W    \>\oplus\>
 U^\prime \otimes W_2) \\ = U_1 \otimes W_2 \>\oplus\> U_2 \otimes W_1.
\end{multline}

To prove (1.4) in the general case, pass to the quotient modulo
$(U_1 \cap U_2) \otimes W + U \otimes (W_1 \cap W_2)$
and obtain by the just proved $U_1 \otimes W_2 + U_2 \otimes W_1$.
\end{proof}

Below, in numerous applications of Lemma 1.1, we will, by abuse of language,
say
\textsl{``by (1.2) and (1.2)$^\prime$, one gets a decomposition 
$I = I_1 \cup I_2 \cup I_3 \cup I_4$ such that 
$S\varphi_i = S^\prime\varphi_i = 0$ for $i \in I_1$, 
$S\varphi_i = T^\prime\alpha_i = 0$ for $i \in I_2$,
$S^\prime\varphi_i = T\alpha_i = 0$ for $i \in I_3$ and 
$T\alpha = T^\prime\alpha_i = 0$ for $i \in I_4$''}.
This means that one can find
a new expression $\Phi = \sum_{i\in I}\varphi_i \otimes \alpha_i$
with indicated properties (where the new $\varphi_i$'s and $\alpha_i$'s
are linear combinations of the old ones).

Unfortunately, for the ``triple intersection'' 
$Ker(S\otimes T) \cap Ker(S^\prime \otimes T^\prime) \cap 
 Ker(S^{\prime\prime} \otimes T^{\prime\prime})$
the analogous decomposition is no longer true. That is why dealing with the second
cohomology group in \S 3, we are unable to obtain a general result and restrict our
considerations with cocycles of rank 1 or with some special cases. For the first cohomology
group, however, Lemma 1.1 suffices to consider the general case, but at the end of
the proof it turns out that it is possible to choose a basis consisting of cocycles of
rank 1.

\section{The first cohomology group}

From now on (in this and subsequent sections), either both $L$ and $M$ or both
$A$ and $V$ are finite-dimensional.

\begin{theorem}\label{2.1}
\begin{multline}\tag{2.1}
H^1(L \otimes A, M \otimes V) \simeq 
H^1(L,M) \otimes V \\ \oplus Hom_L(L,M) \otimes Der(A,V) \oplus
Hom(L/[L,L],M^L) \otimes \frac{Hom(A,V)}{V + Der(A,V)}.
\end{multline}

Each cocycle in $Z^1(L\otimes A, M\otimes V)$ is a linear combination of cocycles
of the three following types (which correspond to the summands in (2.1)):
\begin{enumerate}
\item $x\otimes a \mapsto \varphi(x)\otimes (a\bullet v)$
      for some $\varphi\in Z^1(L,M), v\in V$
\item $x\otimes a \mapsto \varphi(x)\otimes \alpha(a)$
      for some $\varphi\in Hom_L(L,M), \alpha \in Der(A,V)$
\item as in (ii) with $\varphi(L)\subseteq M^L, 
      \varphi([L,L]) = 0, \alpha \in Hom(A,V)$.
\end{enumerate}
\end{theorem}

\begin{remark}
Theorem 2.1 was obtained earlier by Santharoubane \cite{sant} 
in the particular case where $M=L^*, V=A^*$ and $L$ is 1-generated as $U(L)^+$-module,
and by Haddi \cite{haddi-3} (in homological form) in the case of characteristic
zero and $L$ perfect.
\end{remark}

\begin{proof}
Let $\Phi = \sum_{i\in I} \varphi_i \otimes \alpha_i$ be a cocycle of
$$Z^1 (L\otimes A, M\otimes V) \subset Hom(L,M) \otimes Hom(A,V).$$
The cocycle equation $d\Phi = 0$ reads
\begin{equation}\tag{2.2}
\sum_{i\in I} (x\bullet\varphi_i(y) \otimes a\bullet\alpha_i(b) -
               y\bullet\varphi_i(x) \otimes b\bullet\alpha_i(a) -
               \varphi_i([x,y]) \otimes \alpha_i(ab)) = 0.
\end{equation}

Symmetrizing this equation with respect to $x, y$, we get:
$$\notag
\sum_{i\in I} (x\bullet\varphi_i(y) + y\bullet\varphi_i(x)) \otimes
              (a\bullet\alpha_i(b) - b\bullet\alpha_i(a))   = 0.
$$

Substitute $a=b=1$ in (2.2):
$$\notag
\sum_{i\in I} d\varphi_i(x,y) \otimes \alpha_i(1) = 0.
$$

Applying Lemma 1.1 to the last two equations, we get a decomposition 
$I = I_1 \cup I_2 \cup I_3 \cup I_4$ such that
\begin{equation*}
\begin{array}{lll}
d\varphi_i = 0, & x \bullet \varphi_i(y) + y \bullet \varphi_i(x) = 0 & 
\text{for any } i \in I_1,  \\

d\varphi_i = 0, & a \bullet \alpha_i(b) = b \bullet \alpha_i(a) & 
\text{for any } i \in I_2,  \\

x \bullet \varphi_i(y) + y \bullet \varphi_i(x) = 0, & \alpha_i(1) = 0 &
\text{for any } i \in I_3,  \\

a \bullet \alpha_i(b) = b \bullet \alpha_i(a), & \alpha_i(1) = 0 &
\text{for any } i \in I_4. 
\end{array}
\end{equation*}

It is easy to see that $\alpha_i(a) = a \bullet \alpha_i(1)$ for each $i \in I_2$, 
and the mappings 
$x\otimes a \mapsto \varphi_i(x) \otimes \alpha_i(a)$ are cocycles of type (i)
from the statement of the Theorem 2.1, and that $\alpha_i = 0$ for each $i \in I_4$.

Substitute $b=1$ in the cocycle equation (2.2):
$$\notag
\sum_{i\in I_1 \cup I_3} (x\bullet \varphi_i(y) - \varphi_i([x,y])) \otimes
                         (\alpha_i(a) - a\bullet \alpha_i(1)) = 0. 
$$

Now apply Lemma 1.1 again. For elements $\varphi_i$, where $i\in I_1$, the vanishing of
$x\bullet\varphi_i(y) - \varphi_i([x,y])$ implies $\varphi_i([L,L]) = 0$ and 
$\varphi_i(L) \subseteq M^L$, what gives rise to cocycles of type (iii), 
and the vanishing of $\alpha_i(a) - a\bullet\alpha_i(1)$ gives cocycles of type
(i), an already considered case. We have $x\bullet\varphi_i(y) = \varphi_i([x,y])$
for all (remaining) $i \in I_3$.

Hence (2.2) can be rewritten as
$$\notag
\sum_{i\in I_3} \varphi_i([x,y]) \otimes 
                (a\bullet \alpha_i(b) + b\bullet\alpha_i(a) - \alpha_i(ab)) = 0.
$$

The vanishing of the first and second tensor factors gives rise to cocycles of type
(iii) and (ii), respectively.

Hence we have
\begin{multline}\notag
Z^1(L\otimes A, M\otimes V) = Z^1(L,M) \otimes V          \\
+ Hom_L(L,M) \otimes Der(A,V) + Hom(L/[L,L], M^L) \otimes Hom(A,V)
\end{multline}
which can be rewritten as
\begin{multline}\notag
Z^1(L\otimes A, M\otimes V) = Z^1(L,M) \otimes V          \\
\>\oplus\> Hom_L(L,M) \otimes Der(A,V) 
\>\oplus\> Hom(L/[L,L], M^L) \otimes \frac{Hom(A,V)}{V + Der(A,V)}.
\end{multline}

From the considerations above we easily deduce:
$$\notag
B^1 (L\otimes A, M\otimes V) = B^1(L,M) \otimes V
$$
and (2.1) now follows.
\end{proof}

\begin{corollary}\label{2.2}
The derivation algebra of the current Lie algebra $L\otimes A$ is isomorphic to
$$\notag
Der(L)\otimes A \>\oplus\> Hom_L(L,L)\otimes Der(A) \>\oplus\> 
Hom(L/[L,L],Z(L))\otimes \frac{End(A)}{A + Der(A)}.
$$
\end{corollary}

This overlaps with \cite[Theorem 7.1]{block-diff} and \cite[Theorem 1.1]{agg}.

Note that $Hom_L(L,L)$ is nothing but a \textit{centroid} of an algebra $L$
(the set of all linear transformations in $End(L)$ commuting with algebra multiplications).

Specializing to particular cases of $L$ and $A$, 
we get on this way (largely known) results about derivations of some 
particular classes of Lie algebras. So, letting $L = \mathfrak g$, a classical
Lie algebra over $\mathbb C$, and $A = \mathbb C[t,t^{-1}]$, the Laurent polynomial ring,
we get a formula for the derivation algebra of the loop algebra:
$$
Der (\mathfrak g \otimes \mathbb C[t,t^{-1}]) \simeq 
\mathfrak g \otimes \mathbb C[t,t^{-1}] \>\oplus\> 
1 \otimes W,
$$
where $W = Der(\mathbb C[t,t^{-1}])$ is the famous Witt algebra.

More generally, replacing the Laurent polynomial ring by an algebra of functions
meromorphic on a compact Riemann surface and holomorphic outside the fixed
finite set of punctures on the surface, we get a similar formula for the derivation algebra of
a Krichever-Novikov algebra of affine type, where the Witt algebra is replaced by 
a Krichever-Novikov algebra of Witt type.

\section{The second cohomology group}

In this section we obtain some particular results on the second cohomology
group $H^2(L\otimes A, M\otimes V)$. The computations go along the same scheme as for 
$H^1$ but are more complicated.

As we want to express $H^2$ in terms of the tensor products of modules depending on
$(L,M)$ and $(A,V)$, it is natural to do so for underlying modules of the 
Chevalley--Eilenberg complex. We have 
(under the same finiteness assumptions as previous):
\begin{equation}\tag{3.1}\begin{aligned}
C^1(L\otimes A, M\otimes V) &\simeq C^1(L,M) \otimes C^1(A,V)    \\
C^2(L\otimes A, M\otimes V) &\simeq C^2(L,M)\otimes S^2(A,V) \>\oplus\>
                                    S^2(L,M)\otimes C^2(A,V).
\end{aligned}\end{equation}

To obtain a similar decomposition in the third degree, let us denote (by abuse of
language) the Young symmetrizer corresponding to tableau $\lambda$ by the same symbol
$\lambda$. We have decomposition of the unit element in the group algebra $K[S_3]$:
$$\notag
e = 
\frac 16 \>\young(1,2,3)\> + \frac 13 (\>\young(13,2)\> + \>\young(12,3)\>) + \frac 16 \>\young(123)\>.
$$

Then\footnote[2]{
Added April 30, 2010: the arguments at the following few lines (as in the 
published version) are either misleading, or wrong. However, the decomposition 
of $\bigwedge^3(L\otimes A)$ we want to establish is correct and is a particular
case of the Cauchy formula for $n=3$, see \S \ref{sketch}.
I am grateful to Semyon Konstein for pointing this out.
},
using the natural isomorphism $i: T^3(L\otimes A) \simeq T^3(L)\otimes T^3(A)$
and the projection $p: T^3(L\otimes A) \to \wedge^3(L\otimes A)$,
one can decompose the third exterior power of the tensor product as follows:
\begin{multline}\notag
\wedge^3 (L\otimes A) = p \circ (e \times e) \circ i(T^3(L\otimes A)) \\ \simeq
\wedge^3(L) \otimes S^3(A) \oplus 
( \>\young(12,3)\>(L) \otimes \>\young(13,2)\>(A) + 
  \>\young(13,2)\>(L) \otimes \>\young(12,3)\>(A) ) \oplus S^3(L)\otimes \wedge^3(A)
\end{multline}
(all other components appearing in $T^3(L) \otimes T^3(A)$ vanish under the projection).
One directly verifies that
$$\notag
\young(12,3) \times \young(13,2)\> (u) = 0 \quad\text{if and only if}\quad
\young(13,2) \times \young(12,3)\> (u) = 0
$$
for each $u\in \wedge^3(L\otimes A)$.

Hence we get a (noncanonical) isomorphism:
$$\notag
\wedge^3(L\otimes A) \simeq \wedge^3(L) \otimes S^3(A) \>\oplus\>
                            \young(13,2)\>(L) \otimes \young(12,3)\>(A) \>\oplus\>
                            S^3(L) \otimes \wedge^3(A).
$$

Passing to 
$Hom(\>\cdot\>, M\otimes V) \simeq Hom(\>\cdot\>, M) \otimes Hom(\>\cdot\>, V)$, one gets:
\begin{multline}\tag{3.2}
C^3(L\otimes A, M\otimes V) \\ \simeq
C^3(L,M) \otimes S^3(A,V)        \>\oplus\> 
Y^3(L,M) \otimes \widetilde Y^3(A,V) \>\oplus\>
S^3(L,M) \otimes C^3(A,V),
\end{multline}
where $Y^3(L,M) = Hom(\>\young(13,2)\>(L), M), \widetilde Y^3(A,V) = Hom(\>\young(12,3)\>(A), V)$.

According to (3.1)--(3.2) one can decompose $H^2$ as
\begin{equation}\tag{3.3}
H^2 (L\otimes A, M\otimes V) = (H^2)^\prime \oplus (H^2)^{\prime\prime}
\end{equation}
where $(H^2)^\prime$ are the classes of cocycles lying in $C^2(L,M) \otimes S^2(A,V)$
and $(H^2)^{\prime\prime}$
are the classes of cocycles of the form $\Phi + \Psi$, where
$\Phi \in S^2(L,M) \otimes C^2(A,V), \Psi \in C^2(L,M) \otimes S^2(A,V), \Phi \ne 0$.
We will compute $(H^2)^\prime$ and obtain some particular results on
$(H^2)^{\prime\prime}$ (actually $(H^2)^\prime$ and $(H^2)^{\prime\prime}$
are limit terms of a certain spectral sequence; see \S 4).

The differentials of the low degree in the piece 
$$
C^1(L\otimes A, M\otimes V) \overset{d^1}\to C^2(L\otimes A, M\otimes V)
                             \overset{d^2}\to C^3(L\otimes A, M\otimes V)
$$
of the standard Che\-va\-lley-Ei\-len\-berg complex 
can be decomposed as follows:
\begin{align*}
d^1 &= d_1 + d_2       \\
d^2 &= \sum_{\substack{1 \le i \le 2 \\ 1 \le j \le 3}} d_{ij},
\end{align*}
where
\begin{align*}
d_1    &: C^1(L,M) \otimes C^1(A,V) \to C^2(L,M) \otimes S^2(A,V),            \\
d_2    &: C^1(L,M) \otimes C^1(A,V) \to S^2(L,M) \otimes C^2(A,V),            \\
d_{11} &: C^2(L,M) \otimes S^2(A,V) \to C^3(L,M) \otimes S^3(A,V),            \\
d_{12} &: C^2(L,M) \otimes S^2(A,V) \to Y^3(L,M) \otimes \widetilde Y^3(A,V), \\
d_{13} &: C^2(L,M) \otimes S^2(A,V) \to S^3(L,M) \otimes C^3(A,V),            \\
d_{21} &: S^2(L,M) \otimes C^2(A,V) \to C^3(L,M) \otimes S^3(A,V),            \\
d_{22} &: S^2(L,M) \otimes C^2(A,V) \to Y^3(L,M) \otimes \widetilde Y^3(A,V), \\
d_{23} &: S^2(L,M) \otimes C^2(A,V) \to S^3(L,M) \otimes C^3(A,V).
\end{align*}

Direct computations show:
\begin{multline}\notag
d_1 (\varphi\otimes\alpha) (x_1\otimes a_1, x_2\otimes a_2) = \\
    \frac 12 (x_2\bullet\varphi(x_1) - x_1\bullet\varphi(x_2)) \otimes
             (a_1\bullet\alpha(a_2) + a_2\bullet\alpha(a_1)) -
              \varphi([x_1,x_2]) \otimes \alpha (a_1a_2);
\end{multline}
\begin{multline}\notag
d_2 (\varphi\otimes\alpha) (x_1\otimes a_1, x_2\otimes a_2) = \\
    \frac 12 (x_1\bullet\varphi(x_2) + x_2\bullet\varphi(x_1)) \otimes
             (a_2\bullet\alpha(a_1) - a_1\bullet\alpha(a_2));
\end{multline}
\begin{align*}
d_{11} (\varphi\otimes\alpha) (x_1\otimes a_1, x_2\otimes a_2, x_3\otimes & a_3) = \\
    &\frac 13 (\varphi([x_1,x_2],x_3) + \curvearrowright) \otimes
             (\alpha(a_1a_2,a_3)     + \curvearrowright) \\ -
    &\frac 13 (x_1\bullet\varphi(x_2,x_3) + \curvearrowright) \otimes 
             (a_1\bullet\alpha(a_2,a_3)  + \curvearrowright);
\end{align*}
\begin{multline*}
d_{12} (\varphi\otimes\alpha) (x_1\otimes a_1, x_2\otimes a_2, x_3\otimes a_3) = \\
    (2\varphi([x_1,x_2],x_3) + \varphi([x_1,x_3],x_2) - 
      \varphi([x_2,x_3],x_1))\otimes
    (\alpha(a_1a_2,a_3) - \alpha(a_2a_3,a_1)) \\ +
    (- x_1\bullet\varphi(x_2,x_3) + x_2\bullet\varphi(x_1,x_3) + 
       2x_3\bullet\varphi(x_1,x_2)) \otimes
    (a_1\bullet\alpha(a_2,a_3) - a_3\bullet\alpha(a_1,a_2));
\end{multline*}
\begin{multline*}
\shoveright{d_{13} (\varphi\otimes\alpha) (x_1\otimes a_1, x_2\otimes a_2, x_3\otimes a_3) = 0;}
\end{multline*}
\begin{multline*}
d_{22} (\varphi\otimes\alpha) (x_1\otimes a_1, x_2\otimes a_2, x_3\otimes a_3) = \\
    (2\varphi([x_1,x_2],x_3) + \varphi([x_1,x_3],x_2) - 
      \varphi([x_2,x_3],x_1)) \otimes
    (\alpha(a_1a_2,a_3) - \alpha(a_2a_3,a_1)) \\ +
    (- x_1\bullet\varphi(x_2,x_3) + x_2\bullet\varphi(x_1,x_3)) \otimes
    (a_1\bullet\alpha(a_2,a_3) + a_3\bullet\alpha(a_1,a_2) + 
     2a_2\bullet\alpha(a_1,a_3));
\end{multline*}
\begin{multline*}
d_{23} (\varphi\otimes\alpha) (x_1\otimes a_1, x_2\otimes a_2, x_3\otimes a_3) = \\
    \frac 13 (x_1\bullet\varphi(x_2,x_3) + \curvearrowright) \otimes 
             (a_1\bullet\alpha(a_2,a_3) +  \curvearrowright)
\end{multline*}
(the absence of $d_{21}$ in this list is merely a technical matter:
at a relevant stage of computations, it will be convenient to use the entire
differential $d$ rather than $d_{21}$).

Now the reader should be prepared for a bunch of tedious and cumbersome definitions.
We apologize for this, but our excuse is that  all this stuff provides building
blocks for $H^2(L\otimes A, M\otimes V)$ and one can hardly imagine that it may be
defined in a simpler way. Taking a glance at the expressions below, one can believe
that the general formula for $H^n (L\otimes A, M\otimes V)$ hardly exists -- 
if it does, one should give correct $n$-dimensional generalizations of definitions below
(in a few cases this is evident -- like Harrison or cyclic cohomology, but in most cases
it is not).

\begin{definitions}\hfill
\begin{enumerate}
\item
Define $d^{[\>]}, d^\bullet: Hom(L^{\otimes 2}, M) \to Hom(L^{\otimes 3}, M)$
as follows:
\begin{align*}
d^{[\>]} \varphi (x,y,z) &= \varphi([x,y],z) + \curvearrowright \\
d^\bullet \varphi(x,y,z) &= x\bullet\varphi(y,z) + \curvearrowright .
\end{align*}

\item
Define $\wp, D: Hom(A^{\otimes 2}, V) \to Hom(A^{\otimes 3}, V)$ as follows:
\begin{align*}
\wp \alpha(a,b,c) &= \alpha(ab,c) + \curvearrowright \\
D\alpha(a,b,c) &= a\bullet\alpha(b,c) + \curvearrowright.
\end{align*}

\item
$\mathscr B(L,M) = 
\set{\varphi\in C^2(L,M)}{\varphi([x,y],z) + z\bullet\varphi(x,y) = 0; 
                          d^{[\>]}\varphi(x,y,z) = 0}$.
\smallskip

\item
$Q^2(L,M) = \set{d\psi}{\psi\in Hom(L,M); x\bullet\psi(y) = y\bullet\psi(x)}$;

\noindent $H_M^2(L) = (Z^2(L,M^L) + Q^2(L,M))/Q^2(L, M)$.
\smallskip

\item
$\mathscr K(L,M) = \set{\varphi\in C^2(L,M)}
                       {d^{[\>]}\varphi(x,y,z) = 2x\bullet\varphi(y,z)}$;

\noindent $\mathscr J(L,M) = 
\set{\varphi\in C^2(L,M)}
{\varphi(x,y) = 
\psi([x,y]) - \frac 12 x\bullet\psi(y) + \frac 12 y\bullet\psi(x) \text{ for } \psi\in Hom(L,M)}$;

\noindent $\mathscr H(L, M) = (\mathscr K(L,M) + \mathscr J(L,M))/\mathscr J(L,M)$.
\smallskip

\item 
$\mathscr X(L,M) = \set{\varphi\in C^2(L,M)}
{2\varphi([x,y],z) = z\bullet\varphi(x,y); \\ \varphi([x,y],z) = \varphi([z,x],y)}$.
\smallskip

\item
$\mathscr T(L,M) = \set{\varphi\in C^2(L,M)}
{3\varphi([x,y],z) = 2z\bullet\varphi(x,y); \varphi([x,y],z) = \varphi([z,x],y)}$.
\smallskip

\item
$Poor_-(L,M) = \set{\varphi\in C^2(L,M^L)}{\varphi([L,L],L) = 0}$;

\noindent $Poor_+(L,M) = \set{\varphi\in S^2(L,M^L)}{\varphi([L,L],L) = 0}$.
\smallskip

\item
$Sym^2 (L,M) = \set{\varphi\in S^2(L,M)}{x\bullet\varphi(y,z) = y\bullet\varphi(x,z)}$;

\noindent $SB^2 (L,M) = \set{\varphi\in S^2(L,M)}
{\varphi(x,y) = x\bullet\psi(y) + y\bullet\psi(x) \text{ for } \psi \in Hom(L,M)}$;

\noindent $SH^2(L,M) = (Sym^2(L,M) + SB^2(L,M))/SB^2(L,M)$.
\smallskip

\item
Define an action of $L$ on $Hom(L^{\otimes 2}, M)$ via
$$
z \circ \varphi(x, y) = z\bullet\varphi(x,y) + \varphi([x,z],y) + \varphi(x,[y,z]).
$$
$\mathscr S^2(L,M) = \set{\varphi\in S^2(L,M)^L}
{\varphi([x,y],z) + \curvearrowright = 0}$.
\smallskip

\item
$D(A,V) = \set{\beta\in Hom(A,V)}
{\beta(abc) = a\bullet\beta(bc) - bc\bullet\beta(a) + \curvearrowright}$.
\medskip

\item
$HC^1(A,V) = \set{\alpha\in C^2(A,V)}{\wp\alpha = 0}$.
\medskip

\item
$\mathscr C^2(A,V) = \set{\alpha\in C^2(A,V)}
{\alpha(ac,b) - \alpha(bc,a) + a\bullet\alpha(b,c) - b\bullet\alpha(a,c)
 + 2c\bullet\alpha(a,b) = 0}$.
\medskip

\item
$\mathscr P_-(A,V) = \set{\alpha\in C^2(A,V)}
{\alpha(ab,c) = a\bullet\alpha(b,c) + b\bullet\alpha(a,c)}$;

\noindent $\mathscr P_+(A,V) = \set{\alpha\in S^2(A,V)}
{\alpha(ab,c) = a\bullet\alpha(b,c) + b\bullet\alpha(a,c)}$.
\medskip

\item
$\mathscr A(A,V) = \set{\alpha\in S^2(A,V)}{2D\alpha = \wp\alpha}$.
\end{enumerate}
\end{definitions}

The spaces defined in (xi), (xv) are relevant in computation of $Ker\,d_{11}$ (Lemma 3.2),
the spaces defined in (iii)--(viii), (xiv) are relevant in computation of 
$Ker\,d_{11} \cap Ker\,d_{12}$ (see (3.6)), the spaces defined in (ix) are relevant
in computation for the particular case where $L$ is abelian (Proposition 3.5), 
and the spaces defined in (x), (xii)--(xiii) are relevant in computation of 
the relative cohomology group $H^2(L\otimes A; L, M\otimes V)$ (Proposition 3.8).

\begin{remarks}\hfill
\begin{enumerate}
\item 
$d$ (the Chevalley-Eilenberg differential) $= d^{[\>]} + d^\bullet$.

\item
As $B^2(L,M^L) \subseteq Q^2(L,M)$, there is a surjection
$$
H^2(L) \otimes M^L \to H_M^2(L).
$$

\begin{sloppypar}
\item
If $V = K$, then $HC^1(A,V)$ is just the first-order cyclic cohomology $HC^1(A)$.
\end{sloppypar}

\item
The following relations hold:
\begin{gather*}
Poor_-(L,M) \subseteq \mathscr B(L,M) \subseteq Z^2(L,M),    \\
\mathscr B(L,M) \cap Z^2(L,M^L) = Poor_-(L,M),              \\
\mathscr S^2(L, M) \cap S^2(L,M^L)^L = Poor_+(L,M),          \\
\mathscr C^2(A,V) \cap HC^1(A, V) = \mathscr P_-(A,V),      \\ 
\mathscr Z^2(A,V) \cap \mathscr A(A,V) = \mathscr P_+(A,V), \\
Der(A,V) \subseteq D(A,V).
\end{gather*}
\end{enumerate}
\end{remarks}

\begin{proposition}\label{3.1}
\begin{align*}
(H^2)^\prime &\simeq 
           H^2(L,M)        \otimes V
\>\oplus\> H_M^2(L)        \otimes \frac{Hom(A,V)}{V \oplus Der(A,V)}
\>\oplus\> \mathscr H(L,M) \otimes Der(A,V)                                  \\
&\>\oplus\> \mathscr B(L,M) \otimes \frac{Har^2(A,V)}{\mathscr P_+(A,V)}
\>\oplus\> C^2(L,M)^L      \otimes \mathscr P_+(A,V)                         \\
&\>\oplus\> \mathscr X(L,M) \otimes \frac{\mathscr A(A,V)}{\mathscr P_+(A,V)}
\>\oplus\> \mathscr T(L,M) \otimes \frac{D(A,V)}{Der(A,V)}                   \\
&\>\oplus\> Poor_-(L,M)     \otimes \frac{S^2(A,V)}
{Hom(A,V) + D(A,V) + Har^2(A,V) + \mathscr A(A,V)}.
\end{align*}

Each cocycle which lies in $C^2(L,M)\otimes S^2(A,V)$ is a linear combination of cocycles
of the eight following types (which correspond to the respective direct summands
in the isomorphism):
\begin{enumerate}
\item
$x\otimes a \wedge y\otimes b \mapsto \varphi(x,y)\otimes ab\bullet v$, 
where $\varphi\in Z^2(L,M)$ and $v\in V$;

\item
$x\otimes a \wedge y\otimes b \mapsto \varphi(x,y)\otimes \beta(ab)$,
where $\varphi\in Z^2(L,M^L)$ and $\beta\in Hom(A,V)$;

\item
as in (ii) with $\varphi\in \mathscr K(L,M)$ and $\beta\in Der(A,V)$;

\item
$x\otimes a \wedge y\otimes b \mapsto \varphi(x,y)\otimes \alpha(a,b)$,
where $\varphi\in \mathscr B(L,M)$ and $\alpha\in \mathscr Z^2(A,V)$;

\item
as in (iv) with $\varphi\in C^2(L,M)^L$ and $\alpha\in \mathscr P_+(A,V)$;

\item
as in (iv) with $\varphi\in \mathscr X(L,M)$ and $\alpha\in \mathscr A(A,V)$;

\item
$x\otimes a \wedge y\otimes b \mapsto \varphi(x,y)\otimes 
(3a\bullet\beta(b) + 3b\bullet\beta(a) - 2\beta(ab))$, 
where $\varphi\in \mathscr H(L,M)$ and $\beta\in D(A, V)$;

\item
as in (iv) with $\varphi\in Poor_-(L,M)$ and $\alpha\in S^2(A,V)$.
\end{enumerate}
\end{proposition}

\begin{proof}
We have
\begin{equation}\tag{3.4}
(H^2)^\prime = \frac{Ker\,d_{11} \cap Ker\,d_{12}}{Im\,d_1}.
\end{equation}

We compute the relevant spaces in the subsequent series of lemmas.

\begin{lemma}\label{3.2}
\begin{align*}
Ker\,d_{11} &= Z^2(L,M)\otimes V  \\ 
&+ \set{\varphi\in C^2(L,M)}{2d^{[\>]}\varphi + d^\bullet\varphi = 0}\otimes \mathscr A(A,V) \\
&+ \set{\varphi\in C^2(L,M)}{3d^{[\>]}\varphi + 2d^\bullet\varphi = 0}\otimes D(A,V) \\
&+ \set{\varphi\in C^2(L,M)}{d^{[\>]} = d^\bullet\varphi = 0}\otimes S^2(A,V).
\end{align*}
\end{lemma}

\begin{proof}
Substituting $a_1 = a_2 = a_3 = 1$ into the equation $d_{11}\Phi = 0$
(as usual, $\Phi = \sum_{i\in I} \varphi_i\otimes\alpha_i$), one derives the equality
\begin{equation}\tag{3.5}
\sum_{i\in I} d\varphi_i(x_1, x_2, x_3)\otimes\alpha_i(1,1) = 0
\end{equation}
and a decomposition $I = I_1 \cup I_2$ with $d\varphi_i = 0$ for $i\in I_1$
and $\alpha_i(1,1) = 0$ for $i\in I_2$.

Substituting then $a_2 = a_3 = 1$ into the same equation, one gets
$$
\sum_{i\in I} (3d^{[\>]}\varphi_i + 2d^\bullet\varphi_i) \otimes 
              (\alpha_i(1,a_1) - a_1\bullet\alpha_i(1,1)) = 0
$$
and by Lemma 1.1 there is a decomposition 
$I = I_{11} \cup I_{12} \cup I_{21} \cup I_{22}$ with
\begin{equation*}
\begin{array}{lll}
d\varphi_i = 0, & 3d^{[\>]}\varphi_i + 2d^\bullet\varphi_i = 0 & 
\text{for any }i \in I_{11} \\

d\varphi_i = 0, & \alpha_i(1,a) = a\bullet \alpha_i(1,1) 
& \text{for any }i \in I_{12} \\

3d^{[\>]}\varphi_i + 2d^\bullet\varphi_i = 0, & \alpha_i(1,1) = 0 &
\text{for any }i \in I_{21} \\

\alpha_i(1,1) = 0, & \alpha_i(1,a) = a\bullet\alpha_i(1,1) &
\text{for any }i \in I_{22}.
\end{array}
\end{equation*}

Obviously $d\varphi_i = d^\bullet\varphi_i = 0$ for any $i\in I_{11}$, 
so components with $i\in I_{11}$ lie in $Ker\,d_{11}$, 
and $\alpha_i(1,a) = 0$ for any $i\in I_{22}$.

Further, substituting $a_3 = 1$ in our equation, we get
\begin{multline}\notag
\sum_{i\in I} (2d^{[\>]}\varphi_i + d^\bullet\varphi_i) \otimes 
(\alpha_i(a_1,a_2) - 3a_1\bullet\alpha_i(1,a_2) - 3a_2\bullet\alpha_i(1,a_1) \\
 + 2\alpha_i(1,a_1a_2) + 3a_1a_2\bullet\alpha_i(1,1)) = 0.
\end{multline}

In order to apply Lemma 1.1 again, we join the sets $I_{12}$ and $I_{22}$ 
(with the common defining condition $\alpha_i(1, a) = a\bullet \alpha_i(1,1)$) 
and obtain a decomposition 
$I = I_1^\prime \cup I_2^\prime \cup I_3^\prime \cup I_4^\prime$
such that
\[    
\begin{aligned}[t]
2d^{[\>]}\varphi_i + d^\bullet\varphi_i = 0, &\>\> \alpha_i(1,a) = a\bullet \alpha_i(1,1) \\ 
                                             &\>\> \alpha_i(a,b) = ab\bullet \alpha_i(1,1)\\ 
d^{[\>]}\varphi_i = d^\bullet\varphi_i = 0,  &\>\> \alpha_i(1,1) = 0                      \\
3d^{[\>]}\varphi_i + 2d^\bullet\varphi_i = 0, &\>\>
\alpha_i(a,b) = 3a\bullet\alpha_i(1,b) + 3b\bullet\alpha_i(1,a) - 2\alpha_i(1,ab)
\end{aligned}
\hspace{-9pt}
\begin{aligned}[t]
&\text{for any } i \in I_1^\prime     \\
&\text{for any } i \in I_2^\prime     \\
&\text{for any } i \in I_3^\prime     \\
&                                     \\
&\text{for any } i \in I_4^\prime.
\end{aligned}
\]

Note that components $\varphi_i\otimes\alpha_i$ with $i\in I_3^\prime$ are among
those with $i\in I_{11}$ (and lie in $Ker\,d_{11}$).

Now, since the contribution of terms with $i\in I_4^\prime$ to the left side of
(3.5) vanishes, we may apply Lemma 1.1 again, and obtain a decomposition
$I_1^\prime \cup I_2^\prime = 
I_{11}^\prime \cup I_{12}^\prime \cup I_{21}^\prime \cup I_{22}^\prime$
such that
\begin{equation*}
\begin{array}{lll}
2d^{[\>]}\varphi_i + d^\bullet\varphi_i = 0, & \alpha_i(1,a) = 0 
&\text{for any } i\in I_{12}^\prime \\
d\varphi_i = 0,                              & \alpha_i(a,b) = ab\bullet \alpha_i(1,1)
&\text{for any } i\in I_{21}^\prime,
\end{array}
\end{equation*}
and the two remaining types of components do not contribute to the whole picture:
those with indices from $I_{11}^\prime$ satisfy 
$d^{[\>]}\varphi_i = d^\bullet\varphi_i = 0$, the case covered by previous cases,
and those with indices from $I_{22}^\prime$ vanish, as $\alpha_i (a,b) = ab\bullet \alpha_i(1,1) = 0$.
Moreover, the components with indices from $I_{21}^\prime$ lie in $Ker\,d_{11}$.

The remaining part of the equation $d_{11}\Phi = 0$ now reads:
\begin{align*}
\sum_{i\in I_{12}^\prime \cup I_4^\prime} d^{[\>]}\varphi_i (x_1, x_2, x_3)\otimes 
                                         &(\wp\alpha(a_1, a_2, a_3) - 2D\alpha_i(a_2, a_3)  \\
&+ 3a_1a_2\bullet \alpha_i(1, a_3) - a_3\bullet \alpha_i(1, a_1a_2) + \curvearrowright ) = 0.
\end{align*}
Applying Lemma 1.1 again, and noting that the vanishing of the first tensor factor
in each summand above yields the already considered case $d^{[\>]}\varphi_i = d^\bullet\varphi_i = 0$, 
we obtain that the second tensor factor vanishes for all $i\in I_{12}^\prime \cup I_4^\prime$.

Consequently, we obtain two types of components $\varphi_i\otimes \alpha_i$ lying in $Ker\,d_{11}$:
$$
2d^{[\>]}\varphi_i + d^\bullet\varphi_i = 0; \>\>\wp\alpha_i = 2D\alpha_i
$$
and
$$
3d^{[\>]}\varphi_i + 2d^\bullet\varphi_i = 0; \>\> \wp\alpha_i = \frac 32 D\alpha_i; \>\>
\alpha_i \>\text{satisfies the defining condition for}\> i\in I_4^\prime.
$$
The last two conditions imposed on $\alpha_i$ imply $\alpha_i (1, \cdot ) \in D(A,V)$.

Summarizing all this, we obtain the statement of the Lemma.
\end{proof} 

\begin{lemma}\label{3.3}
\begin{multline}\notag
Ker\,d_{12} = C^2(L,M)\otimes V + 
              \set{\varphi\in C^2(L,M)}{x\bullet\varphi(y,z) = z\bullet\varphi(x,y)} \otimes Hom(A,V) \\
+ \set{\varphi\in C^2(L,M)}
      {\varphi([x,y],z) - \varphi([y,z],x) - x\bullet\varphi(y,z) + z\bullet\varphi(x,y) = 0} 
  \otimes \mathscr Z^2(A,V)         \\
+ \set{\varphi\in C^2(L,M)}
      {x\bullet\varphi(y,z) = z\bullet\varphi(x,y); \varphi([x,y],z) = \varphi([y,z],x)}
  \otimes S^2(A,V).
\end{multline}
\end{lemma}

\begin{proof}
Let $\Phi = \sum_{i\in I} \varphi_i\otimes \alpha_i \in Ker\,d_{12}$. 
Substituting $a_2 = 1$ in the equation $d_{12}\Phi = 0$, one gets:
\begin{multline}\notag
\sum_{i\in I} 
(- x_1\bullet\varphi_i(x_2,x_3) + x_2\bullet\varphi(x_1,x_3) + 2x_3\bullet\varphi_i(x_1,x_2)) \\
\otimes (a_1\bullet\alpha_i(1,a_3) - a_3\bullet\alpha_i(1,a_1)) = 0.
\end{multline}
Hence we have a decomposition $I = I_1 \cup I_2$ such that, for $i\in I_1$, 
the first tensor factor in each summand above vanishes, 
and, for $i\in I_2$, the second one vanishes. Elementary
transformations show that
\begin{align*}
x\bullet \varphi_i(y,z) = z\bullet \varphi_i(x,y) &\quad\text{for any } i \in I_1   \\
\alpha_i (1,a) = a\bullet \alpha_i(1,1)           &\quad\text{for any } i \in I_2.  
\end{align*}
Then substituting $a_3 = 1$ into the same initial equation $d_{12}\Phi = 0$, one gets
\begin{multline}\notag
\sum_{i\in I} (2\varphi_i([x_1,x_2],x_3) + \varphi_i([x_1,x_3],x_2) - \varphi_i([x_2,x_3],x_1) \\
               - x_1\bullet \varphi_i(x_2,x_3) + x_2\bullet\varphi_i(x_1,x_3) + 
                2x_3\bullet\varphi_i(x_1,x_2))                                                 \\
              \otimes (\alpha_i(1,a_1a_2) - \alpha_i(a_1,a_2)) = 0.
\end{multline}
Applying Lemma 1.1 and the fact that the vanishing of the first tensor factor here is
equivalent to the condition 
$\varphi_i([x,y],z) - \varphi_i([y,z],x) - x\bullet\varphi_i(y,z) + z\bullet\varphi_i(x,y) = 0$,
we get a decomposition $I = I_{11} \cup I_{12} \cup I_{21} \cup I_{22}$ such that

\medskip
\begin{flushleft}
\begin{tabular}{lll}
$x\bullet\varphi_i(y,z) = z\bullet\varphi_i(x,y)$, &
$\varphi_i([x,y],z) = \varphi_i([y,z],x)$         &
$\text{for any }i\in I_{11}$                       \\

$x\bullet\varphi_i(y,z) = z\bullet\varphi_i(x,y)$, &
$\alpha_i(a,b) = \alpha_i(1,ab)$                  &
$\text{for any }i\in I_{12}$                       \\

\multicolumn{3}{l}
{
$\varphi_i([x,y],z) - \varphi_i([y,z],x) - x\bullet\varphi_i(y,z) + z\bullet\varphi_i(x,y) 
= 0$,
}                                                  \\
& $\alpha_i(1,a) = a\bullet \alpha_i(1,1)$        & 
$\text{for any }i\in I_{21}$                       \\
$\alpha_i(a,b) = ab\bullet \alpha_i(1,1)$          &&
$\text{for any }i\in I_{22}$.      
\end{tabular}
\end{flushleft}
\smallskip

It is easy to see that the components $\varphi_i\otimes \alpha_i$ with indices belonging to 
$I_{11}$, $I_{12}$ and $I_{22}$, already lie in $Ker\,d_{12}$.

The remaining part of the equation $d_{12}\Phi = 0$ becomes:
$$
\sum_{i\in I_{21}} (2\varphi_i([x_1,x_2],x_3) + \varphi_i([x_1,x_3],x_2) - \varphi_i([x_2,x_3],x_1))
                   \otimes \delta\alpha_i(a_1,a_2,a_3) = 0,
$$
where $\delta$ is Harrison(=Hochschild) differential. Thus there is a decomposition 
$I_{21} = I_1^\prime \cup I_2^\prime$, where $\varphi_i$ for $i\in I_1^\prime$ satisfies
the same relations as for $i\in I_{11}$, and $\alpha_i\in \mathscr Z^2(A,V)$
for any $i\in I_2^\prime$.

Putting all these computations together yields the formula desired (the four summands 
there correspond to the defining conditions for $I_{22}$, $I_{12}$, $I_2^\prime$ and $I_{11}$,
respectively; the sum, in general, is not direct).
\end{proof} 

Elementary but tedious transformations of expressions entering in defining conditions
of summands of $Ker\,d_{11}$ and $Ker\,d_{12}$, allow us to write their intersection as
the following direct sum:
\begin{align*}
Ker\,d_{11} \cap Ker\,d_{12} &\simeq 
Z^2(L,M)\otimes V \>\oplus\> Z^2(L,M^L)\otimes \frac{Hom(A,V)}{V \oplus Der(A,V)}  \notag \\
 &\oplus\> \mathscr K(L,M)\otimes Der(A,V) 
\>\oplus\> \mathscr B(L,M)\otimes \frac{\mathscr Z^2(A,V)}{\mathscr P_+(A,V)}     \notag \\
 &\oplus\> C^2(L,M)^L \otimes \mathscr P_+(A,V) 
\>\oplus\> \mathscr X(L,M) \otimes \frac{\mathscr A(A,V)}{\mathscr P_+(A,V)}      \tag{3.6} \\
 &\oplus\> \mathscr T(L,M) \otimes \frac{D(A,V)}{Der(A,V)}                        \notag \\
 \oplus & \> Poor_-(L,M) \otimes \frac{S^2(A,V)}{Hom(A,V) + D(A,V) + \mathscr Z^2(A,V) + \mathscr A(A,V)}.
\end{align*}

According to (3.4), to compute $(H^2)^\prime$, we must consider the equation 
$\Phi = d_1\Psi$,
where $\Phi\in Ker\,d_{11} \cap Ker\,d_{12}$ and $\Psi\in Hom (L\otimes A, M\otimes V)$, 
which is equivalent to elucidation of all possible cohomological dependencies between the 
obtained classes of cocycles.

\begin{lemma}\label{3.4}
Let:

\smallskip

\noindent $\{\varphi_i\}$ be cohomologically independent cocycles in $Z^2(L,M)$,

\noindent $\{\theta_i\}$ be cocycles in $Z^2(L, M^L)$ independent modulo $Q^2(L,M)$, 

\noindent $\{\kappa_i\}$ be elements of $\mathscr K(L,M)$ independent modulo $\mathscr T(L,M)$, 

\noindent $\{\varepsilon_i\}$ be linearly independent cocycles in $\mathscr B(L,M)$, 

\noindent $\{\rho_i\}$ be linearly independent elements in $C^2(L,M)^L$, 

\noindent $\{\chi_i\}$ be linearly independent elements in $\mathscr X(L,M)$, 

\noindent $\{\tau_i\}$ be linearly independent elements in $\mathscr T(L,M)$,

\noindent $\{\xi_i\}$ be linearly independent cocycles in $Poor_-(L,M)$, 

\noindent $\{v_j\}$ be linearly independent elements in $V$,

\noindent $\{\delta_j\}$ be linearly independent derivations in $Der(A,V)$,

\noindent $\{\beta_j\}$ be mappings in $D(A,V)$ independent modulo $Der(A,V)$,

\noindent $\{\gamma_j\}$ be mappings in $Hom(A,V)$ independent both modulo $Der(A,V)$ and 
                                      modulo mappings $a\mapsto a\bullet v$ for all $v\in V$,

\begin{sloppypar}
\noindent $\{F_j\}$ be cocycles in $\mathscr Z^2(A,V)$ independent both cohomologically and 
                                      modulo $\mathscr P_+(A,V)$,
\end{sloppypar}

\noindent $\{P_j\}$ be linearly independent elements in $\mathscr P_+(A,V)$,

\noindent $\{A_j\}$ be elements in $\mathscr A(A,V)$ independent modulo $\mathscr P_+(A,V)$,

\noindent $\{G_j\}$ be mappings in $S^2(A,V)$ independent simultaneously modulo:

\noindent mappings $a\wedge b \mapsto \gamma(ab)$ for all $\gamma\in Hom(A,V)$,

\noindent mappings $a\wedge b \mapsto 3a\bullet \beta(b) + 3b\bullet \beta(a) - 2\beta(ab)$
                for all $\beta\in D(A,V)$, and

\noindent $\mathscr Z^2(A,V) + \mathscr A(A,V)$.

\smallskip

Then the elements of $Ker\,d_{11} \cap Ker\,d_{12}$: 
\begin{multline*}
\varphi_i \otimes (R_{v_j} \circ m), \>
\theta_i \otimes (\gamma_j \circ m), \>
\kappa_i \otimes (\delta_j \circ m), \>
\varepsilon_i \otimes F_j,           \>
\rho_i \otimes P_j,                  \>
\chi_i \otimes A_j,                  \>  \\
\tau_i \otimes (3\delta\beta_j - \beta_j \circ m), \>
\xi_i \otimes G_j
\end{multline*}
($m$ stands for multiplication in A and $R_v$ is an element in $Hom(A,V)$ defined by
$a\mapsto a\bullet v$), are independent modulo $Im\,d_1$.
\end{lemma}


\begin{proof}
We must prove that if
\begin{multline*}
\begin{aligned}
&\sum \varphi_i(x,y)     \otimes ab\bullet v_j    + 
\sum \theta_i (x,y)     \otimes \gamma_j(ab)      +
\sum \kappa_i (x,y)     \otimes \delta_j(ab)   \\ + 
&\sum \varepsilon_i(x,y) \otimes F_j(a,b)         +
\sum \rho_i (x,y)       \otimes P_j(a,b)          +
\sum \chi_i (x,y)       \otimes A_j(a,b)       \\ +
&\sum \tau_i(x,y)        \otimes (3a\bullet\beta_j(b) + 3b\bullet\beta_j(a) - 2\beta_j(ab)) +
\sum \xi_i(x,y)         \otimes G_j(a,b)
\end{aligned}
\\ = 
\sum_{i\in I} \bigg( \psi_i([x,y]) \otimes \alpha_i(ab) + 
\frac 12(-x\bullet\psi_i(y) + y\bullet\psi_i(x)) \otimes (a\bullet\alpha_i(b) + b\bullet\alpha_i(a)) \bigg)
\tag{3.7}
\end{multline*}
for some $\sum_{i\in I} \psi_i\otimes \alpha_i \in Hom(L,M) \otimes Hom(A,V)$
(the right side here is the generic element in $Im\,d_1$), 
then all terms in the left side vanish.

One has $\delta_j(1) = \beta_j(1) = P_j(1,a) = A_j(1,a) = 0$ and one may assume that
$\gamma_j(1) = F_j(1,a) = G_j(1,a) = 0$. Substitute $a = b = 1$ in (3.7):
$$
\sum \varphi_i(x,y) \otimes v_j = \sum_{i\in I} d\psi_i(x,y) \otimes \alpha_i(1).
$$
As $\varphi_i$'s are cohomologically independent and $v_j$'s are linearly independent, the last
equality implies that all summands $\varphi_i(x,y)\otimes v_j$ vanish and there is a decomposition
$I = I_1 \cup I_2$ with $d\psi_i = 0$ for $i\in I_1$ and $\alpha_i(1) = 0$ for $i\in I_2$.

Now substitute $b = 1$ in (3.7):
\begin{multline}\tag{3.8}
\sum \theta_i(x,y) \otimes \gamma_j(a) + 
\sum \kappa_i(x,y) \otimes \delta_j(a) +
\sum \tau_i  (x,y) \otimes \beta_j(a)  \\ =
\sum_{i\in I} (\psi_i([x,y]) + \frac 12(-x\bullet\psi_i(y) + y\bullet\psi_i(x))) \otimes
              (\alpha_i(a) - a\bullet\alpha_i(1)).
\end{multline}
Substituting (3.8) in (3.7), one gets:
\begin{multline*}
  \sum \varepsilon_i(x,y) \otimes F_j(a,b)    + 
  \sum \rho_i       (x,y) \otimes P_j(a,b)    + 
  \sum \chi_i       (x,y) \otimes A_j(a,b) \\ +
3 \sum \tau_i       (x,y) \otimes \delta\beta_j(a,b) + 
  \sum \xi_i        (x,y) \otimes G_j(a,b) \\ =
\frac 12 
\sum (-x\bullet \psi_i(y) + y\bullet \psi_i(x)) \otimes (\delta\alpha_i(a,b) - ab\bullet\alpha_i(1)).
\end{multline*}
The independence conditions of Lemma imply that all summands in the left side
vanish and, due to Lemma 1.1, for $\sum_{i\in I} \psi_i\otimes \beta_i$,
there exists a decomposition $I = I_{11} \cup I_{12} \cup I_{21} \cup I_{22}$ with
\begin{equation*}
\begin{array}{lll}
d\psi_i = 0, &\> x\bullet \psi_i(y) = y\bullet \psi_i(x)     &\text{for any } i\in I_{11}  \\
d\psi_i = 0, &\> \delta\alpha_i(a,b) = ab\bullet \alpha_i(1) &\text{for any } i\in I_{12}  \\
x\bullet \psi_i(y) = y\bullet \psi_i(x), &\> \alpha_i(1) = 0 &\text{for any } i\in I_{21}  \\
\alpha_i(1) = 0, &\> \delta\alpha_i(a,b) = ab\bullet \alpha_i(1) &\text{for any } i\in I_{22}.
\end{array}
\end{equation*}
Denoting $\alpha_i^\prime(a) = \alpha_i(a) - a\bullet \alpha_i(1)$ for $i\in I_{12}$, 
we get $\alpha_i^\prime \in Der(A,V)$.

Substituting all this information back into (3.8), one finally obtains
\begin{multline*}
\sum \theta_i(x,y) \otimes \gamma_j(a) +
\sum \kappa_i(x,y) \otimes \delta_j(a) \\ =
\frac 12 \sum_{i\in I_{12}} \psi_i([x,y]) \otimes (\alpha_i(a) - a\bullet \alpha_i(1)) +
         \sum_{i\in I_{21}} \psi_i([x,y]) \otimes \alpha_i(a) \\ +
\sum_{i\in I_{22}} (\psi_i([x,y]) + \frac 12(-x\bullet\psi_i(y) + y\bullet\psi_i(x))) \otimes \alpha_i(a).
\end{multline*}
The independence conditions of Lemma imply that all terms appearing in the last
equality vanish, and the desired assertion follows.
\end{proof} 

\textit{Conclusion of the proof of Proposition 3.1.}

Lemma 3.3 implies that
\begin{multline*}
Im\,d_1 \simeq B^2(L,M) \otimes V \>\oplus\> 
               (Q^2(L,M) \cap Z^2(L,M^L)) \otimes \frac{Hom(A,V)}{V \oplus Der(A,V)} \\ \>\oplus\>
               (\mathscr T(L,M) \cap \mathscr K(L,M)) \otimes Der(A,V) \>\oplus\>
               \mathscr B(L,M)\otimes Der(A,V)
\end{multline*}
which together with (3.6) entails the asserted isomorphism.
\end{proof} 

Now we turn to computation of the second summand in (3.3), $(H^2)^{\prime\prime}$.

We are unable to compute it in general (and are in doubt about the existence of
a closed general formula for $(H^2)^{\prime\prime}$) and confine ourselves to two
particular cases (in both of them it turns out that $(H^2)^{\prime\prime}$ coincides with the classes
of cocycles lying in $S^2(L,M) \otimes C^2(A,V)$).

\begin{proposition}\label{3.5}
Suppose $L$ is abelian. Then
\begin{align*}
(H^2)^{\prime\prime} \simeq S^2(L, M^L) &\otimes
     \frac {C^2(A,V)}{\set{a\bullet\beta(b) - b\bullet\beta(a)}{\beta\in Hom(A,V)}}  \\
     \oplus\> SH^2(L,M) &\otimes \set{a\bullet\beta(b) - b\bullet\beta(a)}{\beta\in Hom(A,V)}.
\end{align*}
\end{proposition}

First we establish a lemma valid in the general situation (where $L$ is not necessarily abelian).

\begin{lemma}\label{3.6}
\hfill
\begin{enumerate}
\item
$Ker\,d_{23} = Ker\,d^\bullet \otimes C^2(A,V) \\ + 
    S^2(L,M) \otimes \set{\alpha\in C^2(A,V)}{\alpha(a,b) = a\bullet\beta(b) - b\bullet\beta(a)}$
\textnormal{;}
\item 
$Im\,d_2 = \set{\varphi\in S^2(L,M)}{\varphi(x,y) = x\bullet\psi(y) + y\bullet\psi(x)} \\ \otimes
          \set{\alpha\in C^2(A,V)}{\alpha(a,b) = a\bullet\beta(b) - b\bullet\beta(a)}$.
\end{enumerate}
\end{lemma}

\begin{proof}
The only thing which perhaps needs a proof here is the equality
$$
Ker\,D = \set{\alpha\in C^2(A,V)}{\alpha(a,b) = a\bullet\beta(b) - b\bullet\beta(a)}.
$$
The validity of it is verified by appropriate substitution of 1's.
\end{proof} 

\begin{proof}[Proof of Proposition 3.5]
Let $\Phi = \sum_{i\in I} \varphi_i\otimes \alpha_i \in Ker\,d_{23}$, 
with a decomposition on the set of indices $I = I_{1} \cup I_2$ such that
\begin{align*}
d^\bullet \varphi_i (x_1, x_2, x_3) = 0& \quad\text{for any } i\in I_1    \\
\alpha_i(a, b) = -a\bullet\beta_i(b) + b\bullet\beta_i(a)& \quad\text{for any } i\in I_2.
\end{align*}

By Lemma 3.6(i), we may also assume that elements $\alpha_i$, where $i\in I_1$, are independent
modulo $\{a\bullet\beta(b) - b\bullet\beta(a)\}$, and hence $\alpha_i(1,a) = 0$ for each $i\in I_1$.

Suppose there is 
\begin{equation}\tag{3.9}
\Psi = \sum_{i\in I^\prime} \varphi_i^\prime \otimes \alpha_i^\prime \in C^2(L,M) \otimes S^2(A,V)
\end{equation}
such that the class of $\Phi -\Psi$ belongs to $(H^2)^{\prime\prime}$ .
This, in particular, means that $d_{22}\Phi = d_{12}\Psi$:
\begin{multline*}
\sum_{i\in I} (-x_1\bullet\varphi_i(x_2,x_3) + x_2\bullet\varphi_i(x_1,x_3)) \\
   \otimes (a_1\bullet\alpha_i(a_2,a_3) + a_3\bullet\alpha_i(a_1,a_2) + 2a_2\bullet\alpha_i(a_1,a_3)) \\
= \sum_{i\in I^\prime} 
( - x_1\bullet\varphi_i^\prime(x_2,x_3) 
  + x_2\bullet\varphi_i^\prime(x_1,x_3) 
  + 2x_3\bullet\varphi_i^\prime(x_1,x_2)) \\
\otimes (a_1\bullet\alpha_i^\prime(a_2,a_3) - a_3\bullet\alpha_i^\prime(a_1,a_2)).
\end{multline*}
Substituting here $a_2 = 1$, one gets
\begin{multline*}
  2 \sum_{i\in I_1} (-x_1\bullet\varphi_i(x_2,x_3) + x_2\bullet\varphi_i(x_1,x_3)) 
                  \otimes \alpha_i(a_1,a_3)                    \\
+ 3 \sum_{i\in I_1} (-x_1\bullet\varphi_i(x_2,x_3) + x_2\bullet\varphi_i(x_1,x_3))
                  \otimes \alpha_i(a_1,a_3)                    \\
= \sum_{i\in I^\prime} (- x_1\bullet\varphi_i^\prime(x_2,x_3) 
                        + x_2\bullet\varphi_i^\prime(x_1,x_3)
                        + 2x_3\bullet\varphi_i^\prime(x_1,x_2)) \\
              \otimes (a_1\bullet\alpha_i^\prime(1,a_3) - a_3\bullet\alpha_i^\prime(1,a_1)).
\end{multline*}
Hence, due to the independence condition imposed on $\alpha_i$ for $i\in I_1$,
\begin{equation}\tag{3.10}
- x_1\bullet\varphi_i(x_2,x_3) + x_2\bullet\varphi_i(x_1,x_3) = 0,  \quad i\in I_1.
\end{equation}
This, together with condition $\varphi_i\in Ker\,d^\bullet$, evidently implies 
$\varphi_i(L,L) \subseteq M^L$ for
each $i\in I_1$. Note that the terms from $S^2(L,M^L)\otimes C^2(A,V)$ lie in 
$Z^2(L\otimes A, M\otimes V)$.

Now write the cocycle equation for elements from 
$S^2(L,M) \otimes \{a\bullet\beta(b) - b\bullet\beta(a)\}$:
\begin{multline*}
\sum_{i\in I_2} (x_1\bullet\varphi_i(x_2,x_3) 
                \otimes (a_1a_2\bullet\beta_i(a_3) - a_1a_3\bullet\beta_i(a_2))  \\
+ x_2\bullet\varphi_i(x_1,x_3) \otimes (-a_1a_2\bullet\beta_i(a_3) + a_2a_3\bullet\beta_i(a_1)) \\
+ x_3\bullet\varphi_i(x_1,x_2) \otimes (a_1a_3\bullet\beta_i(a_2) - a_2a_3\bullet\beta_i(a_1))) = 0.
\end{multline*}
Substituting $a_2 = a_3 = 1$, we get
$$
\sum_{i\in I_2} (x_2\bullet\varphi_i(x_1,x_3) - x_3\bullet\varphi_i(x_1,x_2)) 
                \otimes (\beta_i(a_1) - a_1\bullet\beta_i(1)) = 0.
$$
As the vanishing of the second tensor factor here leads to the vanishing of the
whole $\alpha_i$, we see that the condition (3.10) holds also in this case, i.e., for all $i\in I_2$.
Conversely, if (3.10) holds, then the cocycle equation is satisfied. Thus the space
of cocycles in $Z^2(L\otimes A, M\otimes V)$ whose cohomology classes lie in $(H^2)^{\prime\prime}$, 
coincides with
\begin{multline*}
S^2(L, M^L) \otimes \frac{C^2(A,V)}{\set{a\bullet\beta(b) - b\bullet\beta(a)}{\beta\in Hom(A,V)}}  \\
\oplus\> \frac{Sym^2(L,M) + SB^2(L,M)}{SB^2(L,M)} \otimes 
\set{a\bullet\beta(b) - b\bullet\beta(a)}{\beta\in Hom(A,V)}
\end{multline*}
(note that we can always take $\Psi = 0$ in (3.9)).

To conclude the proof, one can observe that all these cocycles are cohomologically
independent. This is proved in a pretty standard way, as in Lemma 3.4.
\end{proof} 

Summarizing Proposition 3.1 (for the case where $L$ is abelian) and Proposition 3.5, we obtain

\begin{theorem}\label{3.7}
Let $L$ be an abelian Lie algebra. Then
\begin{align*}
H^2(L\otimes A, M\otimes V) &\simeq H^2(L,M) \otimes V \>\oplus\> \mathscr H(L,M) \otimes Der(A,V) \\
                            &\oplus\> C^2(L,M^L) \otimes \frac{S^2(A,V)}{V \oplus Der(A,V)}  \\
&\oplus\> S^2(L,M^L) \otimes \frac{C^2(A,V)}{\set{a\bullet\beta(b) - b\bullet\beta(a)}{\beta\in Hom(A,V)}} \\
&\oplus\> SH^2(L,M) \otimes {\set{a\bullet\beta(b) - b\bullet\beta(a)}{\beta\in Hom(A,V)}}.
\end{align*}
Each cocycle in $Z^2(L\otimes A, M\otimes V)$ is a linear combination of cocycles of the four
following types (which correspond respectively to the first, the sum of the second
and the third, the fourth and the fifth summands in the isomorphism):
\begin{enumerate}
\item 
$x\otimes a \wedge y\otimes b \mapsto 
\varphi(x,y) \otimes ab\bullet v$ for some $\varphi\in Z^2(L,M)$ and $v\in V$;
\item
$x\otimes a \wedge y\otimes b \mapsto 
\varphi(x,y) \otimes \alpha(a,b)$ for some $\varphi\in C^2(L,M^L)$ and 
$\alpha\in S^2(A,V)$;
\item
as in (ii) with $\varphi\in S^2(L,M^L)$ and $\alpha\in C^2(A,V)$;
\item
$x\otimes a \wedge y\otimes b \mapsto 
\varphi(x,y) \otimes (a\bullet\beta(b) - b\bullet\beta(a))$ for some
$\varphi\in Sym^2(L,M)$ and $\beta\in Hom(A,V)$.
\end{enumerate}
\end{theorem}

\begin{remark}
It is easy to see that if $L$ is abelian, then there is inclusion
$\mathscr H(L,M) \subseteq H^2(L,M)$ ($\mathscr H(L,M)$ consists of classes of cocycles taking values in
$M^L$). Hence, singling out appropriate terms from the first three direct
summands in the isomorphism above, we obtain $\mathscr H(L,M) \otimes S^2(A,V)$
as a direct summand of $H^2(L\otimes A, M\otimes V)$.
\end{remark}
\medskip

Now we want to perform another particular computation of the second cohomology group, 
namely, to compute the relative cohomology $H^2(L\otimes A; L, M\otimes V)$.

We easily see that all constructions can be restricted to the relative complex 
$Hom(\wedge^\star (L\otimes A/K1), M\otimes V)$ with a single 
(but greatly simplifying the matter)
difference that all mappings from $C^3(A,V)$, $Y^3(A,V)$ and $S^3(A,V)$ vanish whenever 
one of their arguments is 1.

We write $(H_L^2)^\prime$ and $(H_L^2)^{\prime\prime}$ to denote the corresponding components
of $H^2(L\otimes A; L, M\otimes V)$.

\begin{proposition}\label{3.8}
\begin{multline*}
(H_L^2)^{\prime\prime} \simeq S^2(L, M^L)^L \otimes HC^1(A,V) \>\oplus\> 
                   \mathscr S^2(L,M) \otimes \frac{\mathscr C^2(L,M)}{\mathscr P_-(A,V)} \\
\oplus\> \frac{S^2(L,M)^L}{S^2(L,M^L)^L} \otimes \mathscr P_-(A,V) \>\oplus\> 
         Poor_+(L,M) \otimes \frac{C^2(A,V)}{HC^1(A,V) + \mathscr C^2(A,V)}
\end{multline*}
\end{proposition}

\begin{proof}
The proof goes along the same scheme as of Proposition 3.1. By Lemma 3.6(i), 
$Ker\,d_{23} = Ker\,d^\bullet \otimes C^2(A,V)$ (as the second tensor factor in the second 
component there vanishes in this case).

The condition $d_{22}\Phi = d_{12}\Psi$ for $\Phi = \sum \varphi_i \otimes \alpha_i \in Ker\,d_{23}$
and $\Psi = \sum \varphi_i^\prime \otimes \alpha_i^\prime \in C^2(L,M) \otimes S^2(L,M)$ reads:
\begin{multline*}
\sum_{i\in I} (2\varphi_i([x_1,x_2],x_3) + \varphi_i([x_1,x_3],x_2) - \varphi_i([x_2,x_3],x_1)) \\
              \otimes (\alpha_i(a_1a_2,a_3) - \alpha_i(a_2a_3,a_1))              \\
+ (- x_1\bullet \varphi_i(x_2,x_3) + x_2\bullet \varphi_i(x_1,x_3)) \otimes
  (a_1\bullet\alpha_i(a_2,a_3) + a_3\bullet\alpha_i(a_1,a_2) + 2a_2\bullet\alpha_i(a_1,a_3)) \\
= \sum_{i\in I^\prime} (2\varphi_i^\prime([x_1,x_2],x_3) + \varphi_i^\prime([x_1,x_3],x_2) -
                         \varphi_i^\prime([x_2,x_3],x_1))  \\
                       \otimes (\alpha_i^\prime(a_1a_2,a_3) - \alpha_i^\prime(a_2a_3,a_1))   \\
+ (-x_1\bullet\varphi_i^\prime(x_2,x_3) + x_2\bullet\varphi_i^\prime(x_1,x_3) + 
   2x_3\bullet\varphi_i^\prime(x_1,x_2))                   \\
\otimes (a_1\bullet\alpha_i^\prime(a_2,a_3) - a_2\bullet\alpha_i^\prime(a_1,a_2)).
\end{multline*}
Substituting here $a_2 = 1$, we obtain (remember about vanishing of all $\alpha$'s if one of
arguments is 1):
\begin{multline*}
\sum_{i\in I} (2\varphi_i([x_1,x_2],x_3) + \varphi_i([x_1,x_3],x_2) - \varphi_i([x_2,x_3],x_1) \\
- x_1\bullet\varphi_i(x_2,x_3) + x_2\bullet\varphi_i(x_1,x_3)) \otimes \alpha_i(a_1,a_3) = 0.
\end{multline*}
This implies
\begin{multline}\tag{3.11}
2\varphi_i([x_1,x_2],x_3) + \varphi_i([x_1,x_3],x_2) - \varphi_i([x_2,x_3],x_1)  \\
    - x_1\bullet\varphi(x_2,x_3) + x_2\bullet\varphi_i(x_1,x_3) = 0, \quad i\in I.
\end{multline}
Since $\varphi\in Ker\,d^\bullet$,
\begin{equation}\tag{3.12}
x_1\bullet\varphi_i(x_2,x_3) + x_3\bullet\varphi_i(x_1,x_3) + x_3\bullet\varphi_i(x_1,x_2) = 0,
\quad i\in I.
\end{equation}
With the help of elementary transformations, (3.11) and (3.12) yield
$$
\varphi_i([x_1,x_3],x_2) + \varphi_i([x_2,x_3],x_1) + x_3\bullet\varphi_i(x_1,x_2) = 0
$$
or, in other words, $\varphi_i \in S^2(L,M)^L$ for each $i\in I$.

Now, writing the cocycle equation for 
$\sum_{i\in I} \varphi_i\otimes\alpha_i \in S^2(L,M)^L \otimes C^2(A,V)$, one
gets
\begin{multline*}
\sum_{i\in I} 
\varphi_i([x_1,x_2],x_3) 
\otimes (\alpha_i(a_1a_2,a_3) - a_1\bullet\alpha_i(a_2,a_3) -a_2\bullet\alpha_i(a_1,a_3))   \\
+ \varphi_i([x_1,x_3],x_2) 
\otimes (-\alpha_i(a_1a_3,a_2) - a_1\bullet\alpha_i(a_2,a_3) + a_3\bullet\alpha_i(a_1,a_2)) \\
+ \varphi_i([x_2,x_3],x_1) 
\otimes (\alpha_i(a_2a_3,a_1) + a_2\bullet\alpha_i(a_1,a_3) + a_3\bullet\alpha_i(a_1,a_2)) = 0.
\end{multline*}

Antisymmetrize this expression with respect to $a_1, a_2$:
\begin{multline*}
\sum_{i\in I} (\varphi_i([x_1,x_3],x_2) + \varphi_i([x_2,x_3],x_1)) 
\otimes (-\alpha_i(a_1a_3,a_2) + \alpha_i(a_2a_3,a_1)                \\
- a_1\bullet\alpha_i(a_2,a_3) + a_2\bullet\alpha_i(a_1,a_3) + 2a_3\bullet\alpha_i(a_1,a_2)) = 0.
\end{multline*}
Consequently, we have a decomposition $I = I_1 \cup I_2$ with
\begin{gather}
\varphi_i([x_1,x_3],x_2) + \varphi_i([x_2,x_3],x_1) = 0, \quad i\in I_1  \tag{3.13}  \\
\alpha_i \in \mathscr C^2(A,V),                          \quad i\in I_2. \tag{3.14}
\end{gather}
Note that (3.13) together with condition $\varphi_i\in S^2(L,M)^L$ implies
$\varphi_i(L,L) \subseteq M^L$ for any $i\in I_1$. Applying to the condition (3.14) the symmetrizer 
$e - (13) + (123)$, we get
\begin{multline}\tag{3.15}
a_1\bullet\alpha_i(a_2,a_3) + a_2\bullet\alpha_i(a_1,a_3) \\
= \frac 13 (2\alpha_i(a_1a_2,a_3) - \alpha_i(a_2a_3,a_1) - \alpha_i(a_1a_3,a_2)), \quad i\in I_2.
\end{multline}
Taking into account (3.13)--(3.15), the cocycle equation can be rewritten as
\begin{multline*}
\sum_{i\in I} (\varphi_i([x_1,x_2],x_3) - \varphi_i([x_1,x_3],x_2) + \varphi_i([x_2,x_3],x_1)) \\
              \otimes (\alpha_i(a_1a_2,a_3) + \alpha_i(a_1a_3,a_2) + \alpha_i(a_2a_3,a_1)) = 0.
\end{multline*}
By Lemma 1.1, there is a decomposition $I = I_{11} \cup I_{12} \cup I_{21} \cup I_{22}$ such that
\begin{equation*}
\begin{array}{lll}
\varphi_i([x,y],z) = \varphi_i(x,[y,z]), & \varphi_i([x,y],z) + \>\curvearrowright\> = 0
&\text{for any } i\in I_{11} \\
\varphi_i([x,y],z) = \varphi_i(x,[y,z]), & \alpha_i \in HC^1(A,V)
&\text{for any } i\in I_{12}  \\
\varphi_i([x,y],z) + \>\curvearrowright\> = 0, & \alpha_i\in \mathscr C^2(A,V)
&\text{for any } i\in I_{21}  \\
\alpha_i\in \mathscr C^2(A,V) \cap HC^1(A,V) &&\text{for any } i \in I_{22}.
\end{array}
\end{equation*}
Evidently, $\varphi_i([L,L],L) = 0$ for any $i\in I_{11}$ and 
$\alpha_i\in \mathscr P_-(A,V)$ for any $i\in I_{22}$. All these four types of components are cocycles in $Z^2(L\otimes A, M\otimes V)$.

Therefore, the space of cocycles whose cohomology classes lie in 
$(H_L^2)^{\prime\prime}$ is as follows:
\begin{multline*}
^L Z^{02} \simeq 
  S^2(L,M^L)^L \otimes HC^1(A,V) 
+ \mathscr S^2(L,M) \otimes \mathscr C^2(A,V) 
+ S^2(L,M)^L \otimes \mathscr P_-(A,V)  \\
+ Poor_+(L,M) \otimes C^2(A,V).
\end{multline*}
(the four summands here correspond to the components indexed by $I_{12}, I_{21}, I_{22}$ and
$I_{11}$ respectively; note that, in this case, we may let $\Psi = 0$ again).

Rewriting this as a direct sum, we get:
\begin{multline*}
S^2(L,M^L)^L \otimes HC^1(A,V) 
\>\oplus\> \mathscr S^2(L,M) \otimes \frac{\mathscr C^2(A,V)}{\mathscr P_-(A,V)}  \\
\oplus\> \frac{S^2(L,M)^L}{S^2(L,M^L)^L} \otimes \mathscr P_-(A,V) 
\>\oplus\> Poor_+(L,M) \otimes \frac{C^2(A,V)}{HC^1(A,V) + \mathscr C^2(A,V)}.
\end{multline*}

And finally, one may show in the same fashion as previously, that all these
cocycles are cohomologically independent, and the assertion of the Proposition follows.
\end{proof} 

Summarizing Propositions 3.1 and 3.8, we obtain:
$$
H^2(L\otimes A; L, M\otimes V) \simeq (H_L^2)^\prime \oplus (H_L^2)^{\prime\prime}
$$
where
\begin{multline*}
(H_L^2)^\prime \simeq \mathscr B(L,M) \otimes \frac{Har^2(A,V)}{\mathscr P_+(A,V)} 
\>\oplus\> C^2(L,M)^L\otimes \mathscr P_+(A,V) 
\>\oplus\> \mathscr X(L,M) \otimes \frac{\mathscr A(A,V)}{\mathscr P_+(A,V)}  \\
\oplus\> Poor_-(L,M) \otimes \frac{S^2(A,V)}{Hom(A,V) + D(A,V) + Har^2(A,V) + \mathscr A(A,V)}
\end{multline*}
and $(H_L^2)^{\prime\prime}$ is described by Proposition 3.8.

We conclude this section with enumeration (for the case of generic $L$) of all possible 
cocycles of rank 1, i.e. those which can be written in the form 
$\varphi\otimes\alpha \in Hom(L^{\otimes 2}, M) \otimes Hom(A^{\otimes 2}, V)$.

In view of (3.3), Propositions 3.1 and 3.8, it suffices to consider cocycles of
rank 1 whose cohomology classes lie in $(H^2)^{\prime\prime}$ and which are independent modulo
$(H_L^2)^{\prime\prime}$. Let us denote this space of cocycles by $Z^{\prime\prime}$.

\begin{proposition}\label{3.9}
Each element of $Z^{\prime\prime}$ is cohomologic to the sum of cocycles of the following
two types:
\begin{enumerate}
\item
$x\otimes a\wedge y\otimes b \mapsto 
\varphi(x,y) \otimes (a\bullet\beta(b) - b\bullet\beta(a))$, where
$\varphi\in Sym^2(L,M)$ is such that $\varphi([L,L],L) = 0$, and $\beta\in Hom(A,V)$;
\item
as in (i) with $\varphi\in S^2(L,M)$, where 
$2\varphi([x,y],z) = x\bullet\varphi(y,z) - y\bullet\varphi(x,z)$, and 
$\beta\in Der(A,V)$.
\end{enumerate}
\end{proposition}

\begin{proof}
Mainly repetition of previous arguments.
\end{proof}

Therefore, there are, in general, 13 types of cohomologically independent cocycles of rank 1 
(7 coming from Proposition 3.1 + 4 coming from Proposition 3.8 + 
2 coming from Proposition 3.9). Of course, in particular cases
some of these types of cocycles may vanish.

We see that, for $H^2(L\otimes A; L, M\otimes V)$ and for 
$H^2(L\otimes A, M\otimes V)$, $L$ abelian, it is possible (in both cases)
to choose a basis consisting of rank 1 cocycles.
In general this is, however, not true. The case of 
$H^2(W_1(n)\otimes A, W_1(n)\otimes A)$, where $W_1(n)$ is the Zassenhaus algebra of positive characteristic, 
treated in \cite{deformations}, shows that there are 
cocycles of rank 2 not cohomologic to (any sum of) cocycles of rank 1.

\section{A sketch of a spectral sequence}\label{sketch}

The computations performed in preceding sections can be described (and
generalized) in terms of a certain spectral sequence. Let us indicate briefly the
main idea (hopefully, the full treatment with further applications will appear elsewhere).

One has the Cauchy formula 
$$
\wedge^n(L\otimes A) \simeq \bigoplus_{\lambda \vdash n} Y_\lambda(L) \otimes Y_{\lambda^\sim}(A),
$$
where $Y_\lambda$ is the Schur functor associated with the Young diagram $\lambda$, and
$\lambda^\sim$ is the Young diagram obtained from $\lambda$ by interchanging
its rows and columns (see, e.g., \cite[p.~121]{fulton-book}).

Applying the functor 
$Hom (\>\cdot\>, M\otimes V) \simeq Hom(\>\cdot\>, M) \otimes Hom(\>\cdot\>, V)$
to both sides of this isomorphism one gets a decomposition
of the underlying modules in the Chevalley--Eilenberg complex:
\begin{equation}\tag{4.1}
C^n(L\otimes A, M\otimes V) \simeq \sum_{\lambda \vdash n} C_\lambda(L,M) \otimes 
C_{\lambda^\sim}(A,V),
\end{equation}
where $C_\lambda(U,W) = Hom(Y_\lambda(U), W)$. The two extreme terms here are 
$C^n(L,M) \otimes S^n(A,V)$ and $S^n(L,M) \otimes C^n(A,V)$.

So each differential $d: C^n(L\otimes A, M\otimes V) \to C^{n+1}(L\otimes A, M\otimes V)$ in the
Chevalley-Eilenberg complex decomposes according to (4.1) into components
$$
d_\lambda^{\lambda^\prime}: C_{\lambda^\prime}(L,M) \otimes C_{\lambda^{\prime\sim}}(A,V) 
\to C_\lambda(L,M) \otimes C_{\lambda^\sim}(A,V)
$$
for each pair $\lambda^\prime \vdash n$ and $\lambda \vdash (n+1)$.
Therefore the following graph of all Young diagrams
\begin{diagram}[width=2.5em,height=4em]
&&&&&& \yng(1) &&&&&&&                                       \\ 
&&&&& \ldTo && \rdTo &&&&&                                      \\ 
&&&& \yng(1,1) &&&& \yng(2) &&&&                                \\ 
&&& \ldTo && \rdTo & & \ldTo(6,2)\ldTo && \rdTo &&&                          \\ 
&& \yng(1,1,1) &&&& \yng(2,1) &&&& \yng(3) &&                    \\ 
& \ldTo && \rdTo\rdTo(4,2) && \ldTo(6,2)\ldTo & \dTo & \rdTo && \ldTo(10,2)\ldTo(6,2)\ldTo(4,2)\ldTo && \rdTo &  \\ 
\yng(1,1,1,1) &&&& \yng(2,1,1) && \yng(2,2) && \yng(3,1) &&&& \yng(4)  \\ 
\dots         &&&& \dots       && \dots     && \dots     &&&& \dots
\end{diagram}

\bigskip\noindent
may be interpreted in the following way: each Young diagram $\lambda$ of size $n$ designates
a module $C_\lambda(L,M) \otimes C_{\lambda^\sim}(A,V)$ and an arrow from $\lambda^\prime$ to 
$\lambda$ represents $d_\lambda^{\lambda^\prime}$.

One can prove that nonzero arrows $d_\lambda^{\lambda^\prime}$ are exactly the following: 
all arrows going
``from right to left'' and those going ``from left to right'' for which either $\lambda^\prime$ 
is a column of height $n$ and $\lambda$ is a diagram of size $n+1$ and of the following shape:
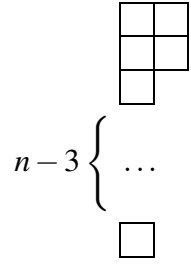
\begin{align*}
                   & \yng(2,2,1)   \\
n-3 \>\Biggl\{\>\> & \dots         \\
                   & \yng(1)
\end{align*}
or $\lambda^\prime$ is included in $\lambda$.

Using this, we can define a decreasing nonnegative filtration $F^kC^\star$ on the complex 
$(C^\star(L\otimes A, M\otimes V), d)$ as the sum of all terms 
$C_\lambda(L,M) \otimes C_{\lambda^\sim}(A,V)$ 
with $\lambda$ belonging to a
``closure'' under nonzero arrows of a single column of height $k+1$.

Now we may consider a (first quadrant) spectral sequence $\{E_r^{\star\star}, d_r \}$ associated
with this filtration. Since the filtration is finite in each degree, the spectral
sequence converges to the desired cohomology group $H^\star(L\otimes A, M\otimes V)$.

Then $E_\infty^{20} = 0$ and $(H^2)^\prime$ and $(H^2)^{\prime\prime}$ from \S 3
are nothing but $E_\infty^{11}$ and $E_\infty^{02}$, respectively.

\section{Structure functions}

In this section we show how the result from \S 3 may be applied to the geometric problem
of calculation of structure functions on manifolds of loops with values in 
compact Hermitian symmetric spaces.

Recall that the base field in this section is $\mathbb C$, what is stipulated by a
geometric nature of the question considered. 
However, all algebraic considerations remain true over any field of 
characteristic $0$.

Let us briefly recall the necessary notions and results. Let $M$ be a complex manifold 
endowed with a $G$-structure (so $G$ is a complex Lie group). \textit{Structure functions}
are sections of certain vector bundles over $M$. 
Their importance stems from the fact that they constitute the complete set 
of obstructions to integrability (= possibility of local flattening) of a given
$G$-structure. In the case $G = O(n)$ structure functions are known under the, perhaps,
more common name \textit{Riemann tensors} (and constitute one of the main objects of study
in the Riemannian geometry). 

A remarkable fact is that structure functions admit a purely algebraic description.
Starting with $\mathfrak g_{-1} = T_m(M)$, the tangent space at a point $m\in M$, and 
$\mathfrak g_0 = Lie(G)$, one may construct, via apparatus of Cartan prolongations, 
a graded Lie algebra $\mathfrak g = \bigoplus_{i\ge -1} \mathfrak g_i$.
Namely, for $i>0$, we have:
\begin{equation}\tag{5.1}
\mathfrak g_i = 
\set{X\in Hom(\mathfrak g_{-1}, \mathfrak g_{i-1})}
{[X(v),w] = [X(w),v] \>\>\text{ for all } v,w \in \mathfrak g_{-1}}.
\end{equation}
 
For any such graded Lie algebra, one may define the \textit{Spencer cohomology groups}
$H_{\mathfrak g_0}^{pq}(\mathfrak g_{-1})$.
Then the space of structure functions of order $k$, i.e. obstructions to identification of the $k$th 
infinitesimal neighborhood of a point $m\in M$ with that of a
point of the manifold with a flat $G$-structure, is isomorphic to the group 
$H_{\mathfrak g_0}^{k2}(\mathfrak g_{-1})$. Note that since
$H^2(\mathfrak g_{-1}, \mathfrak g) = \bigoplus_{k\ge 1} H_{\mathfrak g_0}^{k2}(\mathfrak g_{-1})$,
to compute structure functions for a given $G$-structure on a manifold, one merely needs 
to evaluate the usual Chevalley-Eilenberg cohomology group
$H^2(\mathfrak g_{-1}, \mathfrak g)$ of an abelian Lie algebra $\mathfrak g_{-1}$
with coefficients in the whole $\mathfrak g$ and to
identify structure functions of order $k$ with the graded component 
$\set{\overline{\varphi} \in H^2(\mathfrak g_{-1}, \mathfrak g)}{Im\,\varphi \subseteq \mathfrak g_{k-2}}$,
$k\ge 1$. We refer for the classical text \cite[Chapter VII]{sternberg} for details.

One of the nice examples of manifolds endowed with a $G$-structure are (irreducible) 
compact Hermitian symmetric spaces (CHSS). There are two naturally distinguishable cases:
$rank M = 1$ and $rank M > 1$.

If $rank\,M = 1$, then $M = \mathbb CP^n$, a complex projective space. In this case $\mathfrak g$
turns out to be a general (infinite-dimensional) Lie algebra of Cartan type $W(n)$ 
with a standard grading of depth 1 (recall that $W(n)$ may be defined as a Lie algebra of 
derivations of the polynomial ring in $n$ indeterminates, and consists of
differential operators of the form
$\sum f_i(x_1, \dots, x_n) \partial / \partial x_i$,
$f_i(x_1, \dots, x_n) \in \mathbb C[x_1, \dots, x_n]$).
The result of Serre about cohomology of involutive Lie algebras
of vector fields (see \cite{guillemin-sternberg} for the original Serre's letter and 
\cite{leites-co}, Theorem 1 or \cite{poletaeva-bonn}, p.~9 for a more explicit 
formulation) implies that structure functions in this case vanish.
We will refer for this case as a \textit{rank one case}.

If $rank\,M > 1$, $\mathfrak g$ turns out to
be a classical simple Lie algebra with a grading of depth 1 and length 1:
$\mathfrak g = \mathfrak g_{-1} \oplus \mathfrak g_0 \oplus \mathfrak g_1$.
In particular, Cartan prolongations of order $>1$ vanish, so we might only have structure functions of orders 1, 2 and 3 only
(see \cite{goncharov}, Proposition 4 or \cite{goncharov-selecta}, Proposition 4.2). 
Corresponding structure functions were determined by
Goncharov (\cite{goncharov}, Theorem 1 or \cite{goncharov-selecta}, Theorem 4.5).
We will refer for this case as a \textit{general case}.

\begin{remark}
In the sequel we will need the following well-known fact: for any rank,
\begin{equation}\tag{5.2}
\set{x\in \mathfrak g_i}{[x, \mathfrak g_{-1}] = 0} = 0, \quad i=0,1
\end{equation}
(this condition sometimes is referred as \textit{transitivity} of the corresponding graded Lie 
algebra; see, e.g., \cite{demazure} and references therein).
In particular, $\mathfrak g_{-1}$ is a faithful $\mathfrak g_0$-module.
\end{remark}

During the last decade, there was a big amount of activity 
by Grozman, Leites, Poletaeva, Serganova and Shchepochkina in determining structure functions of various classes of 
(super)manifolds and $G$-structures on them 
(see, e.g., \cite{gls-ium}, \cite{leites-co} and \cite{poletaeva-bonn}
with a transitive closure of references therein).

Here we describe structure functions of manifolds $M^{S^1}$ of
loops with values in a (finite-dimensional) CHSS $M$. The group $G$ here is formally no
longer a Lie group, but its infinite-dimensional analogue, the group of loops,
and the corresponding Lie algebra is a loop Lie algebra 
$\mathfrak g\otimes \mathbb C[t,t^{-1}]$ with a grading inherited from $\mathfrak g$:
$$
\mathfrak g \otimes \mathbb C[t,t^{-1}] = \bigoplus_{i\ge -1} 
\mathfrak g_i\otimes \mathbb C[t,t^{-1}].
$$
The last statement follows from the next simple but handy observation:

\begin{proposition}
Let $\bigoplus_{i\ge -1} \mathfrak g_i$ be the Cartan prolongation of a pair 
$(\mathfrak g_{-1}, \mathfrak g_0)$, where $\mathfrak g_{-1}, \mathfrak g_0$ are finite-dimensional.
Then $\bigoplus_{i\ge -1} (\mathfrak g_i \otimes A)$ is the Cartan prolongation of the pair
$(\mathfrak g_{-1}\otimes A, \mathfrak g_0\otimes A)$.
\end{proposition}

\begin{proof}
Induction on $i$.
As all $\mathfrak g_i$ are finite-dimensional, an element 
$X\in Hom(\mathfrak g_{-1}\otimes A, \mathfrak g_{i-1}\otimes A)$ in the inductive definition (5.1)
of Cartan prolongation may be expressed in the form 
$\sum_{i\in I} \varphi_i \otimes \alpha_i$, where
$\varphi_i\in Hom(\mathfrak g_{-1}, \mathfrak g_{i-1})$, $\alpha_i\in End(A)$.
The rest goes as in the proof of Theorem 2.1.
\end{proof}

Thus, we shall obtain, so to speak, a ``loopization'' of Serre's and Goncharov's results.

In November 1993, Dimitry Leites showed to author a handwritten note by Elena Poletaeva 
containing computations of structure functions of manifolds of loops
corresponding to the following two cases: the (rank one) case $\mathfrak g =  W(1)$
and the (general) case $\mathfrak g = sl(4)$ with graded components 
$\mathfrak g_{-1} = V\otimes V^\star$, 
$\mathfrak g_0 = sl(2) \oplus gl(2), \mathfrak g_1 = V^\star \otimes V$,
where $V$ is a standard 2-dimensional $gl(2)$-module.
Unfortunately, this note has never been published and seems to be lost, 
and more than 10 years 
later nobody from the involved parties cannot recollect the details.
Though formally the main results of this section are generalizations of those Poletaeva's forgotten 
results, it should be noted that 
Poletaeva considered already the typical representatives in both -- 
rank one and general -- cases and observed all the main components and phenomena occurring in
cohomology under consideration.

\vskip 44pt
\begin{definitions}\hfill
\begin{enumerate}
\item
Structure functions (identified with elements of the second cohomology group) generated
by cocycles of the form
$$
(x\otimes a) \wedge (y\otimes b) \mapsto \varphi(x,y) \otimes abu, 
\quad x, y\in \mathfrak g_{-1}, \> a,b \in \mathbb C[t,t^{-1}],
$$
where $\varphi$ is a structure function of CHSS and $u\in \mathbb C[t,t^{-1}]$, will be called \textit{induced}.

\item
Structure functions generated by cocycles of the form
$$
(x\otimes a) \wedge (y\otimes b) \mapsto \varphi(x,y) \otimes \alpha(a,b),
\quad x, y\in \mathfrak g_{-1}, \> a,b \in \mathbb C[t,t^{-1}],
$$
where $\varphi\in C^2(\mathfrak g_{-1}, \mathfrak g_{-1})$ and
$\alpha\in S^2(\mathbb C[t,t^{-1}], \mathbb C[t,t^{-1}])$, will be called 
\textit{almost induced}.

\item
Define a symmetric analogue of $H_{\mathfrak g_0}^{1,2}(\mathfrak g_{-1})$, 
denoted as $SH_{\mathfrak g_0}^{1,2}(\mathfrak g_{-1})$, to be the quotient space
$$
\frac{S^2(\mathfrak g_{-1}, \mathfrak g_{-1})}
{\set{\varphi \in S^2(\mathfrak g_{-1}, \mathfrak g_{-1})}
{\varphi(x,y) = [x,\psi(y)] + [y,\psi(x)] \>\text{for some } \psi\in Hom(\mathfrak g_{-1}, \mathfrak g_0)}}.
$$
\end{enumerate}
\end{definitions}
\bigskip

Clearly, induced and almost induced structure functions arise respectively from
the direct summands 
$H^2(\mathfrak g_{-1}, \mathfrak g) \otimes \mathbb C[t,t^{-1}]$ and 
$\mathscr H(\mathfrak g_{-1}, \mathfrak g) \otimes S^2(\mathbb C[t,t^{-1}], \mathbb C[t,t^{-1}])$ 
of the cohomology group \newline
$H^2(\mathfrak g_{-1} \otimes \mathbb C[t,t^{-1}], \mathfrak g \otimes \mathbb C[t,t^{-1}])$ 
(see Remark after Theorem 3.7 and compare with the paragraph after the proof of Proposition 2.2
in \cite{deformations}).

\begin{theorem}\label{5.2}
For the manifold $M^{S^1}$ of loops with values in a CHSS $M$, the following hold:
\begin{enumerate}
\item
Structure functions can be only of order 1, 2 or 3.
\item
The space of structure functions of order 1 modulo almost induced structure functions is isomorphic to
\begin{multline*}
B_{\mathfrak g_0}^{1,2}(\mathfrak g_{-1}) \otimes 
\frac{S^2(\mathbb C[t,t^{-1}], \mathbb C[t,t^{-1}])}
{(\mathbb C1 \oplus \mathbb C\frac d{dt}) \otimes \mathbb C[t,t^{-1}]} 
\>\oplus\> S^2(\mathfrak g_{-1}, \mathfrak g_{-1}) \otimes 
\frac{C^2(\mathbb C[t,t^{-1}], \mathbb C[t,t^{-1}])}{End\,(\mathbb C[t,t^{-1}])}  \\
\oplus\> SH_{\mathfrak g_0}^{1,2}(\mathfrak g_{-1}) \otimes 
\frac{End(\mathbb C[t,t^{-1}])}{\mathbb C[t,t^{-1}]}.
\end{multline*}
\item
If $rank\,M = 1$, the third direct summand in the last expression vanish.
\item
If $rank\,M = 1$, almost induced structure functions of order 1 and 
all structure functions of order 2 and 3 vanish. 
\item
If $rank\,M > 1$, all structure functions of order 2 and 3 are induced.
\end{enumerate}
\end{theorem}

\begin{remarks}
(i) 
$B_{\mathfrak g_0}^{1,2}(\mathfrak g_{-1})$ is the space of corresponding Spencer 
coboundaries, 
i.e., the space of mappings $\varphi\in C^2(\mathfrak g_{-1}, \mathfrak g_{-1})$ of the form 
$\varphi(x,y) = [x,\psi(y)] - [y,\psi(x)]$ for some $\psi\in Hom\,(\mathfrak g_{-1}, \mathfrak g_0)$.

(ii)
Theorem 3.7 suggests the way in which denominator is embedded into numerator in the three
quotient spaces 
involving $\mathbb C[t,t^{-1}]$ in (ii). In the first quotient space, 
the element 
$(\lambda 1 + \mu\frac d{dt}) t^n \in 
(\mathbb C1 \oplus \mathbb C\frac d{dt}) \otimes \mathbb C[t,t^{-1}]$ 
corresponds to the mapping $\alpha \in S^2(\mathbb C[t,t^{-1}], \mathbb C[t,t^{-1}])$
defined by
$$
\alpha(t^i,t^j) = \lambda t^{i+j+n} + \mu(i + j)t^{i+j+n-1}.
$$
In the second one, the mapping 
$\beta(t^i) = \sum_n \lambda_{in} t^n \in End(\mathbb C[t,t^{-1}])$ 
corresponds to the mapping $\alpha \in C^2(\mathbb C[t,t^{-1}], \mathbb C[t,t^{-1}])$
defined by
$$
\alpha(t^i, t^j) = \sum_n (\lambda_{j,n-i} - \lambda_{i,n-j}) t^n.
$$
In the third one, the element $t^n \in \mathbb C[t,t^{-1}]$ 
corresponds to the mapping $\beta\in End(C[t,t^{-1}])$ which is multiplication by $t^n$:
$$
\beta(t^i) = t^{i+n}.
$$
\end{remarks}

\begin{proof}
Our task is to compute 
$H^2(\mathfrak g_{-1} \otimes \mathbb C[t,t^{-1}], \mathfrak g\otimes \mathbb C[t,t^{-1}])$
for an appropriate $\mathfrak g$.
It turns out that the concrete structure of the Laurent polynomial ring $\mathbb C[t,t^{-1}]$ 
is not important in our approach, and for notational convenience we replace it by an arbitrary
(associative commutative unital) algebra $A$.

Substitute our specific data into the equation of Theorem 3.7:
\begin{align*}
H^2(\mathfrak g_{-1} \otimes A, \mathfrak g \otimes A) 
& \simeq 
H^2(\mathfrak g_{-1}, \mathfrak g) \otimes A 
\>\oplus\> \mathscr H(\mathfrak g_{-1}, \mathfrak g) \otimes Der(A)   \\  \tag{5.3}
& \oplus\> 
C^2(\mathfrak g_{-1}, \mathfrak g^{\mathfrak g_{-1}}) \otimes \frac{S^2(A,A)}{A \oplus Der(A)}  \\
& \oplus\> S^2(\mathfrak g_{-1}, \mathfrak g^{\mathfrak g_{-1}}) \otimes 
\frac{C^2(A,A)}{\set{a\beta(b) - b\beta(a)}{\beta\in End(A)}}         \\
& \oplus\> SH^2(\mathfrak g_{-1}, \mathfrak g) \otimes \set{a\beta(b) - b\beta(a)}{\beta\in End(A)}.
\end{align*}

The next technical lemma is devoted to determination of components appearing
in this isomorphism.

\begin{lemma}
\hfill
\begin{enumerate}
\item $\mathfrak g^{\mathfrak g_{-1}} = \mathfrak g_{-1}$
\item $\mathscr H(\mathfrak g_{-1}, \mathfrak g) = H_{\mathfrak g_0}^{1,2}(\mathfrak g_{-1})$
\item $SH^2(\mathfrak g_{-1}, \mathfrak g) = SH_{\mathfrak g_0}^{1,2}(\mathfrak g_{-1})$.
\end{enumerate}
\end{lemma}

\begin{proof}
(i) Evident in view of (5.2).

(ii) Follows from definitions of appropriate spaces, (5.2) and part (i).

(iii) 
The grading of $\mathfrak g$ induces grading of $SH^2(\mathfrak g_{-1}, \mathfrak g)$:
$$
SH^2(\mathfrak g_{-1}, \mathfrak g) = \bigoplus_{i \ge -1} SH^2_i (\mathfrak g_{-1}, \mathfrak g),
$$
where 
\begin{multline*}
\begin{aligned}
SH^2_i (\mathfrak g_{-1}, \mathfrak g) &=
(Sym^2(\mathfrak g_{-1}, \mathfrak g_i) + SB^2(\mathfrak g_{-1}, \mathfrak g_i)) / 
SB^2(\mathfrak g_{-1}, \mathfrak g_i),     \\
Sym^2(\mathfrak g_{-1}, \mathfrak g_i) &= \set{\varphi\in S^2(\mathfrak g_{-1}, \mathfrak g_i)}
      {[x,\varphi(y,z)] = [y,\varphi(x,z)] \>\text{for all } x,y,z\in \mathfrak g_{-1}},  \\
SB^2(\mathfrak g_{-1}, \mathfrak g_i)  &= 
\{ \varphi\in S^2(\mathfrak g_{-1}, \mathfrak g_i) \>|\> 
   \varphi(x,y) = [x,\psi(y)] + [y,\psi(x)] 
\end{aligned} \\
\text{for some } \psi\in Hom(\mathfrak g_{-1}, \mathfrak g_{i+1}) \}.
\end{multline*}

We immediately see that 
$Sym^2(\mathfrak g_{-1}, \mathfrak g_{-1}) = S^2(\mathfrak g_{-1}, \mathfrak g_{-1})$,
and $\varphi(\cdot,y)$ belongs to the $(i+1)$st Cartan prolongation of the pair
$(\mathfrak g_{-1}, \mathfrak g_0)$
for each $\varphi \in Sym^2(\mathfrak g_{-1}, \mathfrak g_i)$, where $i \ge 0$, and 
$y\in \mathfrak g_{-1}$. 
Hence each $\varphi \in Sym^2(\mathfrak g_{-1}, \mathfrak g_i)$ can be written in the form
\begin{equation}\tag{5.4}
\varphi(x,y) = [x,F(\varphi, y)], \quad \text{for all }x,y \in \mathfrak g_{-1},
\end{equation}
for a certain bilinear map 
$F: Sym^2(\mathfrak g_{-1}, \mathfrak g_i) \times \mathfrak g_{-1} \to \mathfrak g_{i+1}$.

But the symmetry of $\varphi$ implies that 
$F(\varphi, \cdot) \in Hom(\mathfrak g_{-1}, \mathfrak g_{i+1})$ belongs to the
$(i+2)$nd Cartan prolongation of $(\mathfrak g_{-1}, \mathfrak g_0)$.
Hence
\begin{equation}\tag{5.5}
F(\varphi, y) = [y, G(F, \varphi)], 
\quad \text{for all }\varphi\in Sym^2(\mathfrak g_{-1}, \mathfrak g_i), y\in \mathfrak g_{-1},
\end{equation}
for a certain bilinear map 
$G: Hom(Sym^2(\mathfrak g_{-1}, \mathfrak g_i) \times \mathfrak g_{-1}, \mathfrak g_{i+1}) \times 
\mathfrak g_{-1} \to \mathfrak g_{i+2}$.

Combining (5.4) and (5.5) together, one gets $\varphi(x,y) = [x,[y, H(\varphi, y)]]$
for a certain bilinear map 
$H: Sym^2(\mathfrak g_{-1}, \mathfrak g_i) \times \mathfrak g_{-1} \to \mathfrak g_{i+2}$.
Applying again symmetry of $\varphi$, we see that $H$ is constant in the second argument, and hence
each element $\varphi \in Sym^2(\mathfrak g_{-1}, \mathfrak g_i)$ can be written in the form
$\varphi (x,y) = [x,[y,h]]$ for an appropriate $h = H(\varphi, \cdot) \in \mathfrak g_{i+2}$.

But then $\varphi (x,y) = [x, \psi(y)] + [y, \psi(x)]$ for $\psi = - \frac{ad(h)}{2}$, and 
$SH^2_i(\mathfrak g_{-1}, \mathfrak g) = 0$ for $i \ge 0$.

Therefore, $SH^2(\mathfrak g_{-1}, \mathfrak g)$ does not vanish only in the $(-1)$st 
graded component, and the desired equality follows.
\end{proof} 

\textit{Continuation of the proof of Theorem 5.2.}
Substituting the results of
Lemma 5.3 into (5.3), decomposing the Chevalley-Eilenberg cohomology $H^2(\mathfrak g_{-1}, \mathfrak g)$
into the direct sum of corresponding Spencer cohomologies, 
and rearranging the summands as indicated in Remark after Theorem 3.7, we obtain:
\begin{align*}
H^2(\mathfrak g_{-1} \otimes A, \mathfrak g \otimes A) 
& \simeq H_{\mathfrak g_0}^{1,2}(\mathfrak g_{-1}) \otimes S^2(A,A)               \\ 
& \oplus\> (\bigoplus_{k > 1} H_{\mathfrak g_0}^{k2}(\mathfrak g_{-1})) \otimes A \\
& \oplus\> B_{\mathfrak g_0}^{1,2}(\mathfrak g_{-1}) \otimes \frac{S^2(A,A)}{A \oplus Der(A)}  \\
& \oplus\> S^2(\mathfrak g_{-1}, \mathfrak g_{-1}) \otimes 
\frac{C^2(A,A)}{\set{a\beta(b) - b\beta(a)}{\beta\in End(A)}}                     \\
& \oplus\> SH_{\mathfrak g_0}^{1,2}(\mathfrak g_{-1}) \otimes 
{\set{a\beta(b) - b\beta(a)}{\beta\in End(A)}}.
\end{align*}
The first tensor product here consists of almost induced structure functions of order 1 and the second
one consists of induced structure functions of order $>1$, what implies (ii).

As was noted earlier, in the rank one case the first and second tensor product vanish
(what follows from the Serre's theorem), what implies (iv).
In the general case, the second tensor product reduces to structure functions of 
order 2 and 3 -- that is, to
$(H_{\mathfrak g_0}^{2,2}(\mathfrak g_{-1}) \oplus H_{\mathfrak g_0}^{3,2}(\mathfrak g_{-1})) \otimes A$. This proves (i) and (v)
(well, after the final substitution $A = \mathbb C[t,t^{-1}]$).

Part (iii) follows from

\begin{lemma}
For $\mathfrak g = W(n)$ with the standard grading, 
$SH_{\mathfrak g_0}^{1,2}(\mathfrak g_{-1}) = 0$.
\end{lemma}

\begin{proof}
Denoting $\mathfrak g_{-1}$ as $V$, we have $\mathfrak g_0 = gl(V)$, and the statement reduces
to the following: for any $\varphi \in S^2(V,V)$, there is a $\psi \in Hom(V, gl(V))$
such that $\varphi(x,y) = \psi(x)(y) + \psi(y)(x)$. But this is obvious:
take $\psi(x)(y) = \frac 12 \varphi(x,y)$.  
\end{proof}


\begin{remark}
In fact, this trivial reasoning shows that \textit{any} linear mapping $V \times V \to V$,
not necessarily symmetric one, may be represented in the form 
$\varphi(x,y) = \psi(x)(y) + \psi(y)(x)$ for certain $\psi\in Hom(V, gl(V))$.
In particular, it shows that the second Spencer cohomology 
$H_{\mathfrak g_0}^{1,2}(\mathfrak g_{-1})$ vanishes for $\mathfrak g = W(n)$, 
which is a particular case of the
Serre's theorem.
\end{remark}

This completes the proof of the Theorem 5.2.
\end{proof} 

Theorem 5.2 tells how to describe structure functions of manifolds of loops with values in
CHSS in terms of structure functions of underlying CHSS (Spencer cohomology groups), and the space 
$SH_{\mathfrak g_0}^{1,2}(\mathfrak g_{-1})$, which is a sort of a symmetric analogue of 
the Spencer cohomology group.

The thorough treatment of the latter symmetric analogue, including its calculation for 
various $\mathfrak g$'s, 
as well as related construction of a symmetric analogue of Cartan prolongation and some questions
pertained to Jordan algebras and Leibniz cohomology, will, hopefully, appear elsewhere. 
Here we only briefly outline how $SH_{\mathfrak g_0}^{1,2}(\mathfrak g_{-1})$ can be 
determined in the general case (i.e., for classical simple Lie algebras $\mathfrak g$)
in terms of the corresponding root system.

All gradings of length 1 and depth 1 of classical simple Lie algebras may be obtained 
in the following way (see, e.g., \cite{demazure}).
Let $R$ be a root system
of $\mathfrak g$ corresponding to a Cartan subalgebra $\mathfrak h$, $B$ a basis of $R$, 
$\set{h_\beta, e_\alpha}{\beta\in B, \alpha\in R}$ a Chevalley basis of $\mathfrak g$. 
Let $N_{\alpha,\alpha^\prime}$ be structure constants in this basis: 
$[e_\alpha, e_{\alpha^\prime}] = N_{\alpha,\alpha^\prime} e_{\alpha + \alpha^\prime}, 
\alpha + \alpha^\prime \in R$.
Fix a root $\beta\in B$ such that $\beta$
enters in decomposition of each root only with coefficients $-1, 0, 1$ (the existence of
such root implies that $R$ is not of type $G_2$, $F_4$ or $E_8$). Denote by $R_i, i=-1,0,1$,
the set of roots in which $\beta$ enters with coefficient $i$. Then
$$
\mathfrak g_{-1} = \bigoplus_{\alpha\in R_{-1}} \mathbb Ce_\alpha, \quad
\mathfrak g_0    = \mathfrak h \oplus \bigoplus_{\alpha\in R_0} \mathbb Ce_\alpha, \quad
\mathfrak g_1    = \bigoplus_{\alpha\in R_1} \mathbb Ce_\alpha.
$$

Now, consider the mapping
\begin{align*}
T : Hom(\mathfrak g_{-1}, \mathfrak g_0) &\to S^2(\mathfrak g_{-1}, \mathfrak g_{-1}) \\
\psi(x) &\mapsto (T\psi)(x,y) = [x,\psi(y)] + [y,\psi(x)].
\end{align*}
The question of determining $SH_{\mathfrak g_0}^{1,2}(\mathfrak g_{-1})$ 
evidently reduces to evaluation of $Ker\,T$. 

Writing
$$
\psi (e_r) = \sum_{\alpha\in B} \lambda_\alpha^r h_\alpha + \sum_{\alpha\in R_0} \mu_\alpha^r e_\alpha
$$
for $r\in R_{-1}$ and parameters $\lambda_\alpha^r, \mu_\alpha^r \in \mathbb C$, 
we see that the equation
$[x,\psi(y)] + [y,\psi(x)] = 0$ is equivalent to the following three conditions:
\begin{align*}
\sum_{\alpha\in B} \lambda_\alpha^s r(h_\alpha) &= \mu_{r-s}^r N_{s,r-s} 
\quad\text{ for all } r,s\in R_{-1} \text{ such that } r-s \in R_0 \>;     \\
\sum_{\alpha\in B} \lambda_\alpha^s r(h_\alpha) &= 0
\quad\text{ for all } r,s\in R_{-1} \text{ such that } r-s \notin R \>;    \\
\mu_\alpha^r N_{s,\alpha} &= 0 \quad\text{ for all } r,s \in R_{-1}, \alpha\in R_0 \text{ such that } r-s \ne \alpha
\end{align*}
which serve as (linear) defining relations for the space $Ker\,T$ and may be computed in each
particular case.

\section*{Acknowledgements}

This paper was basically written at the beginning of 1990s when I was
a Ph.D. student at Bar-Ilan University under the guidance of Steve Shnider.

My thanks are due to:
Steve Shnider for support and attention;
Dimitry Leites for numerous interesting conversations (in particular, 
I le\-ar\-ned about the structure functions from him) and
useful remarks about preliminary versions of the manuscript (in particular, he spotted a mistake 
in the rank one case in \S 5);
Rutwig Campoamor and Alexander Feldman for facilitating access
to some otherwise hardly available mathematical literature;
Elena Poletaeva for sending her recent reprints;
anonymous referee for very careful reading of the manuscript and numerous remarks
which led to significant improvements (in particular, he spotted a mistake in the rank one 
case in \S 5 as well).

\end{document}